\documentclass[final]{siamltex}
\pdfoutput=1
\usepackage[cp850]{inputenc}
\usepackage[T1]{fontenc}
\usepackage[english]{babel}
\usepackage{graphicx,color}
\usepackage[dvips]{epsfig}
\usepackage{amsmath,amssymb,amsfonts,mathrsfs}
\usepackage{rotating}
\usepackage{cite}
\usepackage[colorlinks=true]{hyperref}

\def\tablenotes{\bgroup\parfillskip=0pt plus 1fil
\leftskip=0pt\relax \rightskip=0pt
\vskip2pt\footnotesize}
\def\endtablenotes{\vskip1pt\egroup}

\def\sphline{\noalign{\vskip3pt}\hline\noalign{\vskip3pt}}

\def\Xint#1{\mathchoice
   {\XXint\displaystyle\textstyle{#1}}%
   {\XXint\textstyle\scriptstyle{#1}}%
   {\XXint\scriptstyle\scriptscriptstyle{#1}}%
   {\XXint\scriptscriptstyle\scriptscriptstyle{#1}}%
   \!\int}
\def\XXint#1#2#3{{\setbox0=\hbox{$#1{#2#3}{\int}$}
     \vcenter{\hbox{$#2#3$}}\kern-.5\wd0}}

\def\dashint{\Xint-}

\newtheorem{algorithm}[theorem]{Algorithm}

\title{On the use of conformal maps for the acceleration of convergence of the trapezoidal rule and Sinc numerical methods}

\author{Richard Mikael Slevinsky\thanks{Department of Mathematical and Statistical Sciences, University of Alberta, Edmonton, Alberta, Canada. (rms8@ualberta.ca)} \and Sheehan Olver\thanks{School of Mathematics and Statistics, The University of Sydney, Sydney, Australia. (Sheehan.Olver@sydney.edu.au)}}

\begin{document}

\maketitle

\begin{abstract}
We investigate the use of conformal maps for the acceleration of convergence of the trapezoidal rule and Sinc numerical methods. The conformal map is a polynomial adjustment to the $\sinh$ map, and allows the treatment of a finite number of singularities in the complex plane. In the case where locations are unknown, the so-called Sinc-Pad\'e approximants are used to provide approximate results. This adaptive method is shown to have almost the same convergence properties. We use the conformal maps to generate high accuracy solutions to several challenging integrals, nonlinear waves, and multidimensional integrals.
\end{abstract}

\begin{keywords}
Trapezoidal rule. Sinc numerical methods. Conformal maps.
\end{keywords}

\begin{AMS}
30C30, 41A30, 65D30, 65L10.%
\end{AMS}

\section{Introduction}

The trapezoidal rule is one of the most well-known methods in numerical integration. While the composite rule has geometric convergence for periodic functions, in other cases it has been used as the starting point of effective methods, such as Richardson extrapolation~\cite{Richardson-210-459-11} and Romberg integration~\cite{Romberg-28-30-55}. The geometric convergence breaks down with endpoint singularities, and this issue inspired a different approach to improve on the composite rule. From the Euler-Maclaurin summation formula, it was noted that some form of exponential convergence can be obtained for integrands which vanish at the endpoints, suggesting that undergoing a variable transformation may well induce this convergence~\cite{Schwartz-4-19-69,Takahasi-Mori-1-201-71,Stenger-12-103-73}. After this observation, the race was on to determine exactly which variable transformation, and therefore which decay rate, is optimal. Numerical experiments showed the exceptional promise of rules such as the $\tanh$ substitution~\cite{Evans-Forbes-Hyslop-15-339-84}, the ${\rm erf}$ substitution~\cite{Takahasi-Mori-21-206-73}, the IMT rule~\cite{Iri-Moriguti-Takasawa-17-3-87}, and the $\tanh$-$\sinh$ substitution~\cite{Takahasi-Mori-9-712-74}, among others~\cite{Mori-14-713-78}. But exactly which one is optimal, and in which setting?

Using a functional analysis approach, this question was beautifully answered by establishing the optimality of a double exponential endpoint decay rate for the trapezoidal rule on the real line for approximating analytic integrands~\cite{Sugihara-75-379-97}. The domain of analyticity is described in terms of a strip of maximal width $\pi$ centred on the real axis in the complex plane. This optimality also prescribed the optimal step size and a near-linear convergence rate ${\cal O}(e^{-kN/\log N})$, where $N$ is the number of sample points and $k$ is a constant proportional to the strip width.

The results allowed for displays of strong performance for integrals with integrable endpoint singularities without changing the rule in any way~\cite{Mori-12-119-85,Mori-41-897-05,Mori-Sugihara-127-287-01,Sugihara-Matsuo-164-673-04}. The double exponential transformation was also adapted to Fourier and general oscillatory integrals in~\cite{Ooura-Mori-38-353-91,Ooura-Mori-112-229-99}. Recognizing the trapezoidal rule as the integration of a Sinc expansion of the integrand, the double exponential advocates adapted their analysis to Sinc approximations~\cite{Sugihara-72-767-03}, and also to all the numerical methods therewith derived, such as Sinc-Galerkin and Sinc-collocation methods~\cite{Horiuchi-Sugihara-5-99,Sugihara-149-239-02,Muhammad-Nurmuhammad-Mori-Sugihara-182-32-05} for initial and boundary value problems, Sinc indefinite integration~\cite{Tanaka-Sugihara-Murota-74-655-04}, iterated integration~\cite{Muhammad-Sugihara-22-77-05}, and Sinc-collocation for integral equations~\cite{Muhammad-Nurmuhammad-Mori-Sugihara-177-269-05,Okayama-Matsuo-Sugihara-234-1211-10}, all obtaining the near-linear convergence rate ${\cal O}(e^{-kN/\log N})$. More recently, the researchers then focused on improving the original convergence results by developing more precise upper estimates on the error given bounds on the function~\cite{Tanaka-Sugihara-Murota-Mori-111-631-09,Tanaka-Sugihara-Murota-78-1553-09,Okayama-Matsuo-Sugihara-124-361-13,Okayama-Tanaka-Matsuo-Sugihara-125-511-13,Tanaka-Okayama-Matsuo-Sugihara-125-545-13}.

In this work, we report an improvement of the trapezoidal rule in the context of a finite number of singularities -- of any kind -- near the contour of integration. This problem has been considered before, in Gaussian quadrature and in Sinc quadrature~\cite{Smith-Lyness-23-231-69,Monegato-47-301-86,Bialecki-29-464-89,Gautschi-64-417-13}. The prevailing philosophy seems to be to characterize singularities as specifically as possible, then account for them by either adding terms from Cauchy's residue theorem to the approximation or by modifying the weights and abscissas. While the examples and applications show exceptional performance of the algorithms, the case of general singularities is still untreated. In the optimal double exponential framework, singularities near the integration contour may reduce the width of the strip of analyticity about the real axis. As the double exponential decay rate is typically induced by a variable transformation, we seek to find variable transformations which place the threatening singularities on the upper and lower edges of the maximal strip of width $\pi$. In this work, such variable transformations are conformal maps which maximize the convergence rate despite the presence of the singularities.

The idea of using conformal maps to speed up numerical computations is not new. In fact, it was recently pioneered by Tee and Trefethen~\cite{Tee-Trefethen-28-1798-06}, Hale and Trefethen~\cite{Hale-Trefethen-46-930-08}, Hale and Tee~\cite{Hale-Tee-31-3195-09}, and Hale's so-titled Ph.~D.~thesis~\cite{Hale-Thesis-09}. Inspired by this work, we investigate the use of the Schwarz-Christoffel map from the strip of width $\pi$ to a polygonally bounded region, with the possibility of some sides being infinite. In this way, the edges of the strip of width $\pi$ contain the pre-images of the function's singularities, their possible branches, and other objects limiting analyticity. In the case dealt with in~\cite{Hale-Tee-31-3195-09}, an algorithm is constructed to solve the Schwarz-Christoffel parameter problem, and the integral definition of the Schwarz-Christoffel map is simplified by partial fraction decomposition.

For the strip map variation of the Schwarz-Christoffel map, the explicit integration of the Schwarz-Christoffel map could not be done~\cite{Hale-Tee-31-3195-09}, so in this work, an approximate map is constructed based on polynomial adjustments to the $\sinh$ map. We choose the $\sinh$ map because it appears in every double exponential map of the canonical finite, infinite, and semi-infinite domains. The polynomial adjustments add exactly enough parameters to locate a finite number of pre-images of singularities on the edges of the maximal strip of width $\pi$, and the parameter problem -- the determination of such polynomials adjustments -- is in complete analogy with the Schwarz-Christoffel parameter problem. However, the approximate map is significantly less expensive to evaluate, and therefore well-paired with the double exponential transformation for high precision numerical experiments.

The problem of poles limiting analyticity has been noted in~\cite{Sugihara-Matsuo-164-673-04}, where Sugihara and Matsuo show the double exponential Sinc expansion of the function:
\begin{equation}\label{eq:SEandDEEx5}
f(x) = \dfrac{x\,(1-x)\,e^{-x}}{(1/2)^2+(x-1/2)^2},\qquad x\in[0,1],
\end{equation}
is not as efficient as methods of polynomial interpolation. To demonstrate how simple our nonlinear program can be, we note that while the original double exponential transformation for the problem is $\phi(t) = \frac{1}{2}\tanh(\frac{\pi}{2}\sinh t)+\frac{1}{2}$, the optimized map is $\phi(t) = \frac{1}{2}\tanh(\frac{\pi}{4}\sinh t)+\frac{1}{2}$, and the convergence rate is approximately tripled.

We demonstrate the merits of our algorithm on four integrals, each with its own combination of singularities. In these cases, the algorithm obtains approximately $2.5$--$4$ times as many correct digits as a na\"ive double exponential transformation for the same number of function evaluations. The algorithm is applied to obtain solutions of the forced Benjamin-Ono equation describing nonlinear waves. Then, the algorithm also shows its merit in the evaluation of multidimensional integrals.

High accuracy in scientific computing is important in many applications. As an example, high accuracy results are required to confidently obtain results from an integer relation algorithm such as PSLQ~\cite{Ferguson-Bailey-92} for finding closed forms for integrals. As well, it is also important to have a rapidly convergent algorithm so that problems are solved in acceptable computational times.

\section{Quadrature and Sinc methods by variable transformation}

Using a variable transformation to induce exponential decay at the endpoints is first performed in~\cite{Takahasi-Mori-21-206-73} and it is extended to double exponential decay in~\cite{Takahasi-Mori-9-712-74}. They find that variable transformations that induce double exponential decay at the endpoints perform better than single exponential transformations.

For infinite and semi-infinite integrals with or without pre-existing exponential decay, various other transformations have been proposed to induce double exponential decay. Examples from~\cite{Tanaka-Sugihara-Murota-Mori-111-631-09} are included in Table~\ref{table:SEandDEmaps}.

\begin{table}[htbp]
\begin{center}
\caption{Variable transformations $\phi(t)$ for endpoint decay.}
\label{table:SEandDEmaps}
\begin{tabular*}{\hsize}{@{\extracolsep{\fill}}lll}
\sphline
Interval & Single Exponential & Double Exponential\\
\sphline
$[-1,1]$ & $\tanh(t/2)$ & $\tanh(\frac{\pi}{2}\sinh t)$\\
$(-\infty,+\infty)$ & $\sinh(t)$ & $\sinh(\frac{\pi}{2}\sinh t)$\\
$[0,+\infty)$ & $\log(\exp(t)+1)$ & $\log(\exp(\frac{\pi}{2}\sinh t)+1)$\\
$[0,+\infty)$ & $\exp(t)$ & $\exp(\frac{\pi}{2}\sinh t)$\\
\sphline
\end{tabular*}
\end{center}
\end{table}
We follow closely the rigorous derivations in~\cite{Sugihara-75-379-97,Sugihara-72-767-03} of the optimality of the trapezoidal rule and Sinc numerical methods for functions with double exponential decay as $x\to\pm\infty$. Let $d$ be a positive number and let $\mathscr{D}_d$ denote the strip region of width $2d$ about the real axis:
\begin{equation}
\mathscr{D}_d = \{z\in\mathbb{C}:|{\rm Im}\,z|<d\}.
\end{equation}
Let $\omega(z)$ be a non-vanishing function defined on the region $\mathscr{D}_d$, and define the Hardy space $H^\infty(\mathscr{D}_d,\omega)$ by:
\begin{equation}
H^\infty(\mathscr{D}_d,\omega) = \{f:\mathscr{D}_d\to\mathbb{C}|\,f(z)\textrm{ is analytic in }\mathscr{D}_d\textrm{, and }||f||<+\infty\},
\end{equation}
where the norm of $f$ is given by:
\begin{equation}
||f|| = \sup_{z\in\mathscr{D}_d}\left|\dfrac{f(z)}{\omega(z)}\right|.
\end{equation}
Consequentially, for $\omega(z)$ that decays double exponentially, the functions in $H^\infty(\mathscr{D}_d,\omega)$ decay double exponentially as well.
Let us consider the $N(=2n+1)$-point trapezoidal rule for the interval $(-\infty,+\infty)$:
\begin{equation}
\int_{-\infty}^{+\infty}f(x){\rm\,d}x \approx h\sum_{k=-n}^{+n}f(k\,h),
\end{equation}
where the mesh size $h$ is suitably chosen for a given positive integer $n$. For the trapezoidal rule, let $\mathscr{E}_{N,h}^{\rm T}(H^\infty(\mathscr{D}_d,\omega))$ denote the error norm in $H^\infty(\mathscr{D}_d,\omega)$:
\begin{equation}
\mathscr{E}_{N,h}^{\rm T}(H^\infty(\mathscr{D}_d,\omega)) = \sup_{||f||\le1}\left|\int_{-\infty}^{+\infty}f(x){\rm\,d}x - h\sum_{k=-n}^{+n}f(k\,h)\right|.
\end{equation}
Let $B(\mathscr{D}_d)$, originally introduced in~\cite{Stenger-12-103-73}, denote the family of all functions $f$ analytic in $\mathscr{D}_d$ such that:
\begin{equation}
\mathcal{N}_1(f,\mathscr{D}_d) = \int_{\partial\mathscr{D}_d}|f(z)|{\rm\,d}z <+\infty.
\end{equation}

\begin{theorem}[Sugihara~\cite{Sugihara-75-379-97}]\label{thm:SEconvergence}
Suppose that the function $\omega(z)$ satisfies the following three conditions:
\begin{enumerate}
\item $\omega(z) \in B(\mathscr{D}_d)$;
\item $\omega(z)$ does not vanish at any point in $\mathscr{D}_d$ and takes real values on the real axis;
\item the decay rate of $\omega(z)$ on the real axis is specified by:
\begin{equation}
\alpha_1\exp\left(-(\beta|x|^\rho)\right)\le|\omega(x)|\le \alpha_2\exp\left(-(\beta|x|^\rho)\right),\quad x\in\mathbb{R},
\end{equation}
where $\alpha_1,\alpha_2,\beta>0$ and $\rho\ge1$.
\end{enumerate}
Then:
\begin{equation}
\mathscr{E}_{N,h}^{\rm T}(H^\infty(\mathscr{D}_d,\omega)) \le C_{d,\omega} \exp\left(-(\pi d\beta N)^{\frac{\rho}{\rho+1}}\right),
\end{equation}
where $N=2n+1$, the mesh size $h$ is chosen optimally as:
\begin{equation}
h = (2\pi d)^{\frac{1}{\rho+1}}(\beta n)^{-\frac{\rho}{\rho+1}},
\end{equation}
and $C_{d,\omega}$ is a constant depending on $d$ and $\omega$.
\end{theorem}

\begin{theorem}[Sugihara~\cite{Sugihara-75-379-97}]\label{thm:DEconvergence}
Suppose that the function $\omega(z)$ satisfies the following three conditions:
\begin{enumerate}
\item $\omega(z) \in B(\mathscr{D}_d)$;
\item $\omega(z)$ does not vanish at any point in $\mathscr{D}_d$ and takes real values on the real axis;
\item the decay rate of $\omega(z)$ on the real axis is specified by:
\begin{equation}
\alpha_1\exp(-\beta_1e^{\gamma|x|})\le|\omega(x)|\le \alpha_2\exp(-\beta_2e^{\gamma|x|}),\quad x\in\mathbb{R},
\end{equation}
where $\alpha_1,\alpha_2,\beta_1,\beta_2,\gamma>0$.
\end{enumerate}
Then:
\begin{equation}
\mathscr{E}_{N,h}^{\rm T}(H^\infty(\mathscr{D}_d,\omega)) \le C_{d,\omega} \exp\left(-\dfrac{\pi d\gamma N}{\log(\pi d\gamma N/\beta_2)}\right),
\end{equation}
where $N=2n+1$, the mesh size $h$ is chosen optimally as:
\begin{equation}
h = \dfrac{\log(2\pi d\gamma n/\beta_2)}{\gamma n},
\end{equation}
and $C_{d,\omega}$ is a constant depending on $d$ and $\omega$.
\end{theorem}

Since the trapezoidal rule is equivalent to the integration of the Sinc expansion of a function~\cite{Stenger-23-165-81}, the entire process of analyzing the convergence rates with different endpoint decay can also be useful for the Sinc expansion of a function, with subtle differences that arise in the mesh size and convergence rates. Let us consider the $N(=2n+1)$-point Sinc approximation of a function on the real line:
\begin{equation}
f(x) \approx \sum_{j=-n}^{+n}f(j\,h)S(j,h)(x),
\end{equation}
where $S(j,h)(x)$ is the so-called Sinc function:
\begin{equation}
S(j,h)(x) = \dfrac{\sin[\pi(x/h-j)]}{\pi(x/h-j)},
\end{equation}
and where the step size $h$ is suitably chosen for a given positive integer $n$. From l'H\^opital's rule, it can easily be seen that the Sinc functions are mutually orthogonal at the so-called Sinc points $x_k = k\,h$:
\begin{equation}
S(j,h)(k\,h) = \delta_{k,j},
\end{equation}
where $\delta_{k,j}$ is the Kronecker delta~\cite{Abramowitz-Stegun-65}.

For the Sinc approximation, let $\mathscr{E}_{N,h}^{\rm Sinc}(H^\infty(\mathscr{D}_d,\omega))$ denote the error norm in $H^\infty(\mathscr{D}_d,\omega)$:
\begin{equation}
\mathscr{E}_{N,h}^{\rm Sinc}(H^\infty(\mathscr{D}_d,\omega)) = \sup_{||f||\le1}\left\{\sup_{x\in\mathbb{R}}\left|f(x) - \sum_{j=-n}^{+n}f(j\,h)\,S(j,h)(x)\right|\right\}.
\end{equation}

\begin{theorem}[Sugihara~\cite{Sugihara-72-767-03}]\label{thm:SESincconvergence}
Suppose that the function $\omega(z)$ satisfies the following three conditions:
\begin{enumerate}
\item $\omega(z) \in B(\mathscr{D}_d)$;
\item $\omega(z)$ does not vanish at any point in $\mathscr{D}_d$ and takes real values on the real axis;
\item the decay rate of $\omega(z)$ on the real axis is specified by:
\begin{equation}
\alpha_1\exp\left(-(\beta|x|^\rho)\right)\le|\omega(x)|\le \alpha_2\exp\left(-(\beta|x|^\rho)\right),\quad x\in\mathbb{R},
\end{equation}
where $\alpha_1,\alpha_2,\beta>0$ and $\rho\ge1$.
\end{enumerate}
Then:
\begin{equation}
\mathscr{E}_{N,h}^{\rm Sinc}(H^\infty(\mathscr{D}_d,\omega)) \le C_{d,\omega} N^{\frac{1}{\rho+1}}\exp\left(-\left(\dfrac{\pi d\beta N}{2}\right)^{\frac{\rho}{\rho+1}}\right),
\end{equation}
where $N=2n+1$, the mesh size $h$ is chosen optimally as:
\begin{equation}
h = (\pi d)^{\frac{1}{\rho+1}}(\beta n)^{-\frac{\rho}{\rho+1}},
\end{equation}
and $C_{d,\omega}$ is a constant depending on $d$ and $\omega$.
\end{theorem}

\begin{theorem}[Sugihara~\cite{Sugihara-72-767-03}]\label{thm:DESincconvergence}
Suppose that the function $\omega(z)$ satisfies the following three conditions:
\begin{enumerate}
\item $\omega(z) \in B(\mathscr{D}_d)$;
\item $\omega(z)$ does not vanish at any point in $\mathscr{D}_d$ and takes real values on the real axis;
\item the decay rate of $\omega(z)$ on the real axis is specified by:
\begin{equation}
\alpha_1\exp(-\beta_1e^{\gamma|x|})\le|\omega(x)|\le \alpha_2\exp(-\beta_2e^{\gamma|x|}),\quad x\in\mathbb{R},
\end{equation}
where $\alpha_1,\alpha_2,\beta_1,\beta_2,\gamma>0$.
\end{enumerate}
Then:
\begin{equation}
\mathscr{E}_{N,h}^{\rm Sinc}(H^\infty(\mathscr{D}_d,\omega)) \le C_{d,\omega} \exp\left(-\dfrac{\pi d\gamma N}{2 \log(\pi d\gamma N/(2\beta_2))}\right),
\end{equation}
where $N=2n+1$, the mesh size $h$ is chosen optimally as:
\begin{equation}
h = \dfrac{\log(\pi d\gamma n/\beta_2)}{\gamma n},
\end{equation}
and $C_{d,\omega}$ is a constant depending on $d$ and $\omega$.
\end{theorem}

Last but not least, there is the nonexistence theorem, which provides a fundamental bound for the proposed optimization approach.
\begin{theorem}[Sugihara~\cite{Sugihara-75-379-97}]\label{thm:DEnonexistence}
There exists no function $\omega(z)$ that satisfies the following three conditions:
\begin{enumerate}
\item $\omega(z)\in B(\mathscr{D}_d)$;
\item $\omega(z)$ does not vanish at any point in $\mathscr{D}_d$ and takes real values on the real axis;
\item the decay rate on the real axis of $\omega(z)$ is specified as:
\begin{equation}
\omega(x) = {\cal O}\left(\exp(-\beta e^{\gamma|x|})\right),\quad {\rm as}\quad |x|\to\infty,
\end{equation}
where $\beta>0$, and $d\gamma>\pi/2$.
\end{enumerate}
\end{theorem}

\section{Maximizing the convergence rates}

From the previous theorems~\ref{thm:DEconvergence} and \ref{thm:DESincconvergence} on the convergence rates of the trapezoidal rule with a prescribed decay at the endpoints and the nonexistence theorem~\ref{thm:DEnonexistence} of analytic functions with double exponential decay in too wide a strip, we may ask the following question. How can we use a conformal map $\phi$ to maximize the convergence rate of the trapezoidal rule:
\begin{equation}
\int_{-\infty}^\infty f(\phi(t))\phi'(t){\rm\,d}t \approx h\sum_{k=-n}^{+n} f(\phi(k\,h))\phi'(k\,h),
\end{equation}
or the Sinc approximation:
\begin{equation}
f(x) \approx \sum_{j=-n}^{+n} f(\phi(j\,h)) S(j,h)(\phi^{-1}(x)),
\end{equation}
despite the singularities of $f\in\mathbb{C}$ which limit its domain of analyticity? To formulate this problem mathematically, let $\Phi_{\rm ad}$ be the admissible space of all functions $\phi$ satisfying the conditions of theorems~\ref{thm:DEconvergence} and~\ref{thm:DESincconvergence}:
$$\Phi_{\rm ad} = \left\{\begin{array}{lll}
\phi & : & \hspace*{-0.2cm}f(\phi(\cdot))\phi'(\cdot)\in H^\infty(\mathscr{D}_d,\omega)\textrm{~for~some~}d>0,\\
&&\textrm{and~for~some~}\omega\textrm{~such~that:}\\
& 1. & \omega(z)\in B(\mathscr{D}_d);\\
& 2. & \omega(z)\textrm{~does~not~vanish~at~any~point~in~}\mathscr{D}_d\\
&& \textrm{~and~takes~real~values~on~the~real~axis};\\
& 3. & \alpha_1\exp\left(-\beta_1e^{\gamma|x|}\right)\le|\omega(x)|\le \alpha_2\exp\left(-\beta_2e^{\gamma|x|}\right),\\
&& x\in\mathbb{R},\textrm{~~where~}\alpha_1,\alpha_2,\beta_1,\beta_2,\gamma>0.
\end{array}\right\}$$

We wish to find the $\phi\in\Phi_{\rm ad}$ so that the convergence rates are maximized:
\begin{align*}
\underbrace{\underset{\phi\in\Phi_{\rm ad}}{\rm argmax} \left(\dfrac{\pi d\gamma N}{\log(\pi d\gamma N/\beta_2)}\right)}_{\rm Convergence~Theorem~\ref{thm:DEconvergence}} &\qquad{\rm subject~to}\quad \underbrace{d\gamma \le \dfrac{\pi}{2}}_{\rm Nonexistence~Theorem~\ref{thm:DEnonexistence}}\\~&\\
\underbrace{\underset{\phi\in\Phi_{\rm ad}}{\rm argmax} \left(\dfrac{\pi d\gamma N}{2\log(\pi d\gamma N/(2\beta_2))}\right)}_{\rm Convergence~Theorem~\ref{thm:DESincconvergence}} &\qquad{\rm subject~to}\quad \underbrace{d\gamma \le \dfrac{\pi}{2}}_{\rm Nonexistence~Theorem~\ref{thm:DEnonexistence}}
\end{align*}
As infinite-dimensional optimization problems for $\phi$, these are challenging problems. However, the convergence rates of theorems~\ref{thm:DEconvergence} and~\ref{thm:DESincconvergence} are asymptotic ones and therefore it is of equivalent interest to investigate the asymptotic solutions to the problem. Consider the asymptotic problems:
\begin{align}
\dfrac{\pi d\gamma N}{\log(\pi d\gamma N/\beta_2)} & = \dfrac{\pi d\gamma N}{\log N + \log(\pi d\gamma/\beta_2)},\nonumber\\
& \sim \dfrac{\pi d\gamma N}{\log N},\quad{\rm as}\quad N\to\infty,\\
\dfrac{\pi d\gamma N}{2\log(\pi d\gamma N/(2\beta_2))} & = \dfrac{\pi d\gamma N}{2\log N + 2\log(\pi d\gamma/(2\beta_2))},\nonumber\\
& \sim \dfrac{\pi d\gamma N}{2\log N},\quad{\rm as}\quad N\to\infty.
\end{align}
Then, the linear appearance of $d\gamma$ leads directly to the following result.
\begin{theorem}\label{thm:DEmaxcon}
Let $\Phi_{\rm as,ad} = \left\{\Phi_{\rm ad}:d\gamma=\pi/2\right\}$ be the asymptotically admissible subspace of the admissible space $\Phi_{\rm ad}$. Then for every $\phi_{\rm as}\in\Phi_{\rm as,ad}$:
\begin{equation}
\mathscr{E}_{N,h}^{\rm T}(H^\infty(\mathscr{D}_d,\omega)) \le C_{d,\omega} \exp\left(-\dfrac{\pi^2 N}{2\log(\pi^2N/2\beta_2)}\right),
\end{equation}
where $N=2n+1$, the mesh size $h$ is chosen optimally as:
\begin{equation}
h = \dfrac{\log(\pi^2 n/\beta_2)}{\gamma n},
\end{equation}
and $C_{d,\omega}$ is a constant depending on $d$ and $\omega$. This same $\phi_{\rm as}$ ensures that:
\begin{equation}
\mathscr{E}_{N,h}^{\rm Sinc}(H^\infty(\mathscr{D}_d,\omega)) \le C_{d,\omega} \exp\left(-\dfrac{\pi^2 N}{4\log(\pi^2N/4\beta_2)}\right),
\end{equation}
where $N=2n+1$, the mesh size $h$ is chosen optimally as:
\begin{equation}
h = \dfrac{\log(\pi^2 n/2\beta_2)}{\gamma n},
\end{equation}
and $C_{d,\omega}$ is a constant depending on $d$ and $\omega$.
\end{theorem}

The implication of such a theorem is that suitable mappings $\phi$ can be found which maximize the convergence rates by neutering the terrible effects of singularities near the approximation interval. In this section, we find such mappings by starting with the observation that in all of the maps in Table~\ref{table:SEandDEmaps} for the finite, semi-infinite, and infinite canonical domains, an elementary map is composed with the $\sinh$ map. Therefore, it seems as though studying the $\sinh$ map, or some modification thereof, will be the best place to start.

Let $f$ have a finite number of singularities located at the points $\{\delta_k\pm{\rm i}\epsilon_k\}_{k=1}^n$. The four maps in Table~\ref{table:SEandDEmaps} can be written as the composition of:
\begin{align}
\psi(z) = \tanh(z), \quad& \psi^{-1}(z) = \tanh^{-1}(z),\label{eq:outerDEmap1}\\
\psi(z) = \sinh(z), \quad& \psi^{-1}(z) = \sinh^{-1}(z),\label{eq:outerDEmap2}\\
\psi(z) = \log(e^z+1), \quad& \psi^{-1}(z) = \log(e^z-1),\label{eq:outerDEmap3}\\
\psi(z) = \exp(z), \quad& \psi^{-1}(z) = \log(z).\label{eq:outerDEmap4}
\end{align}
with the $\frac{\pi}{2}\sinh$ function. In any of these cases, let a finite number of singularities of $f$ be transformed as $\{\tilde{\delta}_k\pm{\rm i}\tilde{\epsilon}_k\}_{k=1}^n$ as the ordered set of $\{\psi^{-1}(\delta_k\pm{\rm i}\epsilon_k)\}_{k=1}^n$, where $\tilde{\delta}_1<\tilde{\delta}_2<\cdots<\tilde{\delta}_n$.

The $\sinh$ function is a conformal map from the strip $\mathscr{D}_{\frac{\pi}{2}}$ to the entire complex plane with two branch cuts emanating outward from the points $\pm{\rm i}$. It is actually the most rudimentary Schwarz-Christoffel formula mapping from the strip $\mathscr{D}_{\frac{\pi}{2}}$ to the entire complex plane~\cite{Howell-Trefethen-11-928-90} with those two aforementioned branches. Let $g$ map the strip $\mathscr{D}_{\frac{\pi}{2}}$ to the polygonally bounded region $P$ having vertices $\{w_k\}_{k=1}^{n} = \{\tilde{\delta}_1+{\rm i}\tilde{\epsilon}_1$, \ldots, $\tilde{\delta}_n+{\rm i}\tilde{\epsilon}_n\}$ and interior angles $\{\pi\alpha_k\}_{k=1}^{n}$. Let also $\frac{\pi}{2}\alpha_{\pm}$ be the divergence angles at the left and right ends of the strip $\mathscr{D}_{\frac{\pi}{2}}$. Then the function:
\begin{equation}
g(z) = A + C\int^z e^{(\alpha_--\alpha_+)\zeta} \prod_{k=1}^n \left[\sinh(\zeta-z_k)\right]^{\alpha_k-1}{\rm\,d}\zeta,
\end{equation}
where $z_k = g(w_k)$ and for some $A$ and $C$ maps the interior of the top half of the strip $\mathscr{D}_{\frac{\pi}{2}}$ to the interior of the polygon $P$. The solution of the constants $A$, $C$, and $\{z_k\}_{k=1}^n$ is known as the Schwarz-Christoffel parameter problem.

\begin{figure}[htbp]
\begin{center}
\begin{tabular}{ccc}
\includegraphics[width=0.3\textwidth]{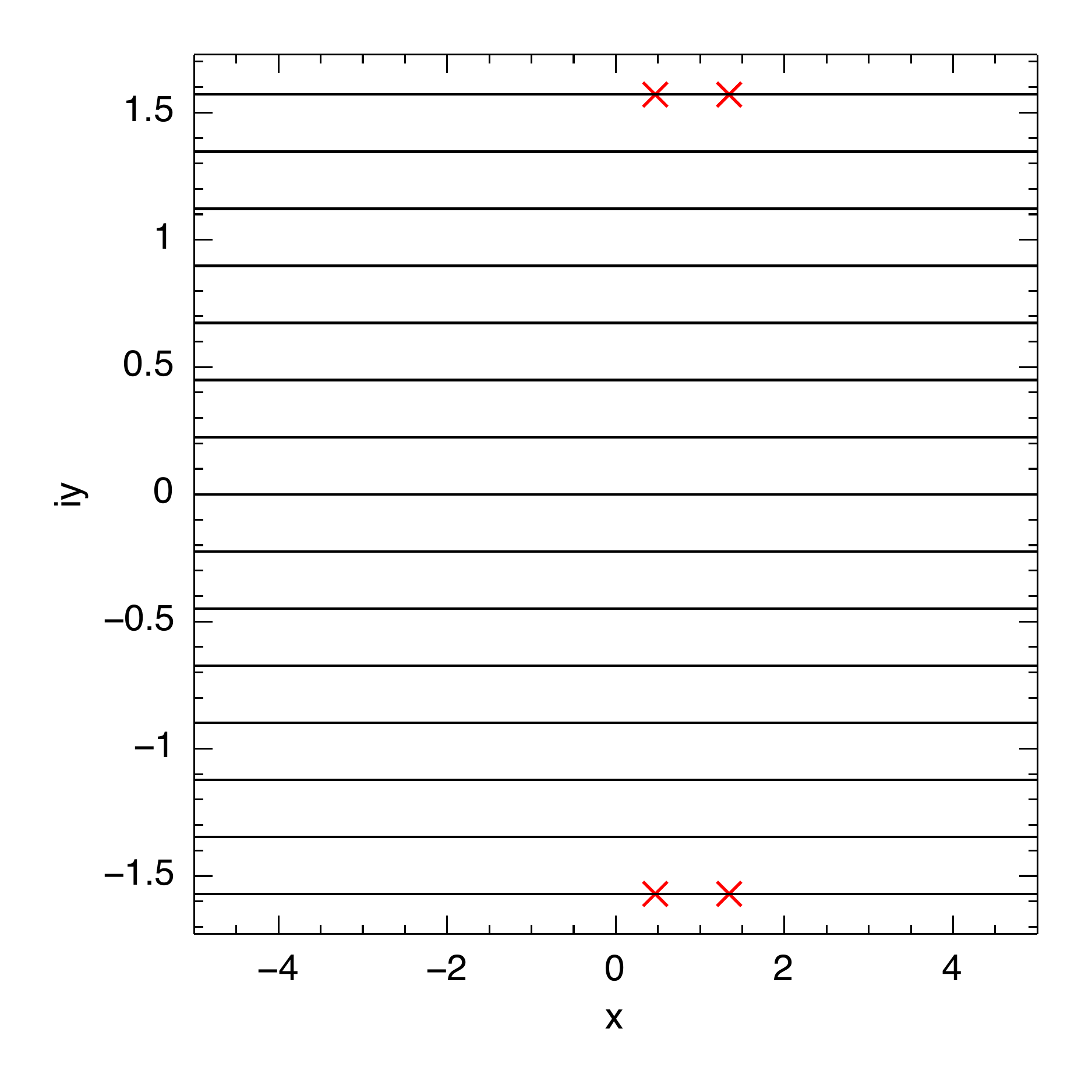}&
\includegraphics[width=0.3\textwidth]{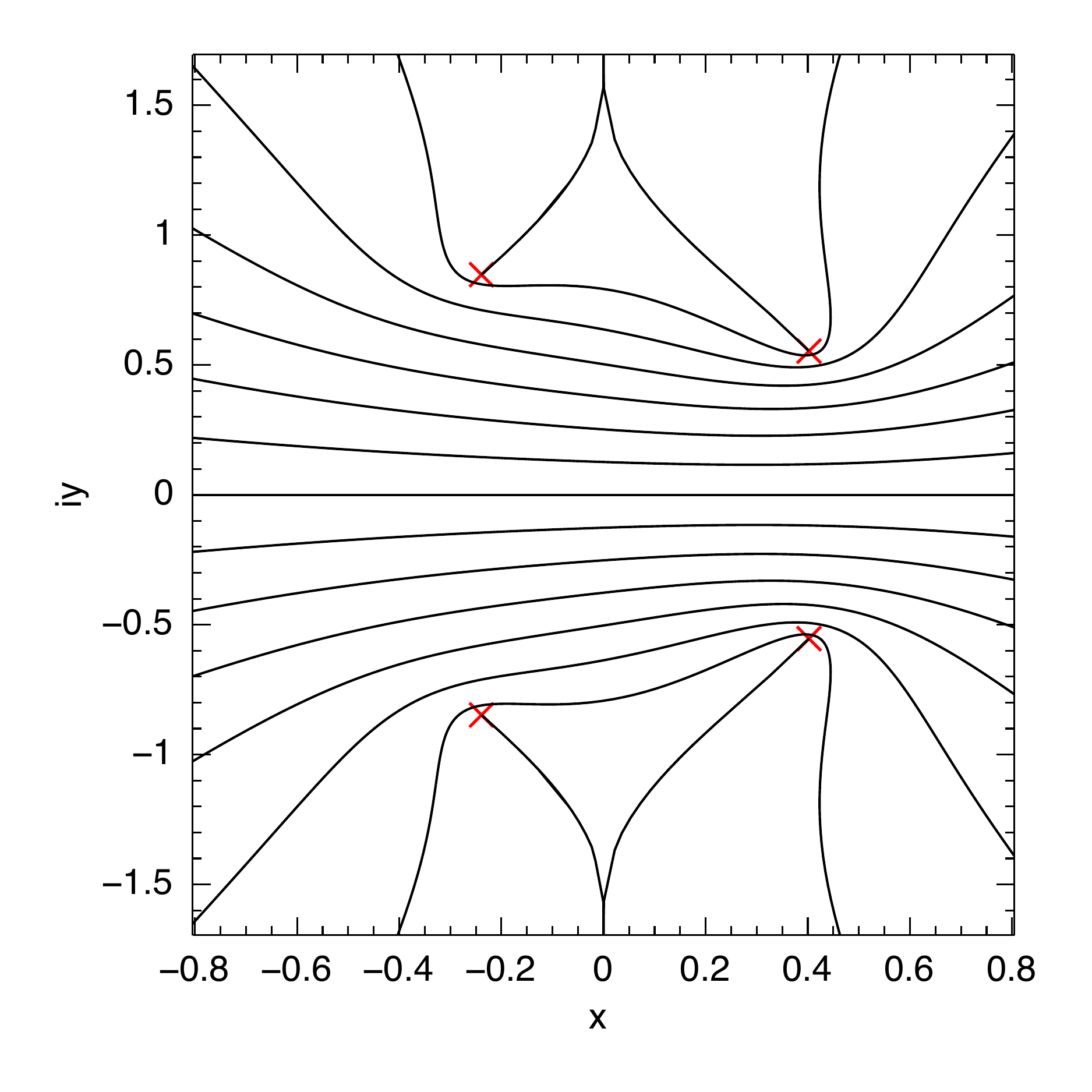}&
\includegraphics[width=0.3\textwidth]{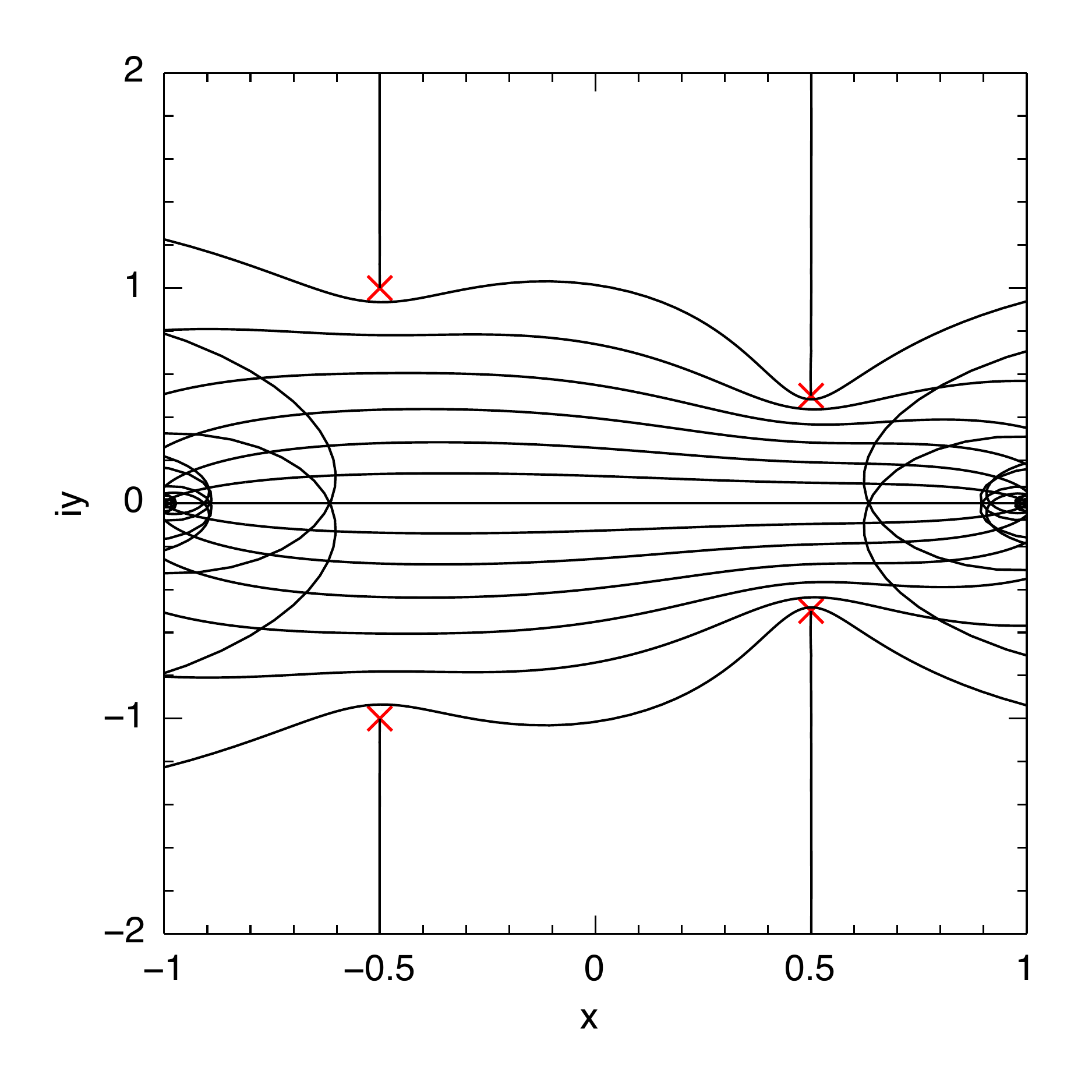}\\
(a) & (b) & (c)\\
\end{tabular}
\caption{In (a) the plot of the strip $\mathscr{D}_{\frac{\pi}{2}}$ with singularities located on the boundary, in (b) the resulting Schwarz-Christoffel map, and in (c) a $\tanh$ map of the Schwarz-Christoffel map. In all three cases, the crosses track the singularities $\delta_1\pm{\rm i}\epsilon_1 = -1/2 \pm {\rm i}$ and $\delta_2\pm{\rm i}\epsilon_2 = 1/2\pm{\rm i}/2$. For the sake of comparison, an integral with these singularities is treated in Example~\ref{example:SEandDEEx1}.}
\label{fig:SCmap}
\end{center}
\end{figure}

Figure~\ref{fig:SCmap} shows an example of the Schwarz-Christoffel map for the polygonal restrictions on $\mathbb{C}$ due to the possible singularities of $f$. The Schwarz-Christoffel map is actually the exact solution of the problem of maximizing the convergence rates, as it maps points on the top and bottom of the strip $\mathscr{D}_{\frac{\pi}{2}}$ to the singularities. However, the entire process is computationally intensive. Firstly, the nonlinear system of equations of the Schwarz-Christoffel parameter problem needs to be solved, and secondly, the map is defined as an integral. The parameter problem can be prohibitive to solve, requiring thousands of integrations of the map function. Also, the integral only has an analytical expression for a polygon with one finite vertex, and this gives the $\sinh$ map. The Schwarz-Christoffel Toolbox in MATLAB~\cite{Trefethen-1-82-80,Trefethen-Driscoll-3-533-98,Driscoll-Trefethen-02} is used to solve for the maps in Figure~\ref{fig:SCmap}, and provides a precision of approximately $10^{-8}$ for a computation time on the order of one minute. In Figure~\ref{fig:SCmap} and subsequent figures, the plots show the mapping of lines with constant imaginary values between $-{\rm i}\pi/2$ and $+{\rm i}\pi/2$ via the conformal map from the strip, then the composition of this conformal map with one of the maps $\psi(z)$ of~\eqref{eq:outerDEmap1}--\eqref{eq:outerDEmap4}. Were it only for the difficulties posed by the Schwarz-Christoffel parameter problem, this approach may have some promise. However, the major problem is that even after the parameter problem is solved, the map itself is defined as an integral and requires a large computational effort compared to the following proposed approach.

Fortunately, due to the framework of the double exponential transformation, we can make a polynomial adjustment to the $\sinh$ map while still retaining a variable transformation $\phi$ which induces double exponential decay. For any real values of the $n+1$ parameters $\{u_k\}_{k=0}^n$, the function:
\begin{equation}\label{eq:SCalternative}
h(t) = u_0\sinh(t) + \sum_{j=1}^{n} u_j t^{j-1},\qquad u_0>0,
\end{equation}
still grows single exponentially. Therefore, the composition $\psi(h(t))$ for any $\psi$ in~\eqref{eq:outerDEmap1}--\eqref{eq:outerDEmap4} still induces a double exponential variable transformation. The benefit of choosing such functions is that we now have sufficient parameters which we can use to ensure the pre-images of the singularities $\{\tilde{\delta}_k\pm{\rm i}\tilde{\epsilon}_k\}_{k=1}^n$ reside on the top and bottom edges of the strip $\mathscr{D}_{\frac{\pi}{2}}$, respectively. This is done by solving the system of equations:
\begin{equation}\label{eq:SCalternativesystem}
h(x_k+{\rm i}\pi/2) = \tilde{\delta}_k + {\rm i}\tilde{\epsilon}_k,\quad{\rm for}\quad k=1,\ldots,n.
\end{equation}
This is a system of $n$ complex equations for the $2n+1$ unknowns $\{u_k\}_{k=0}^n$ and the $x$-coordinates of the pre-images of the singularities $\{x_k\}_{k=1}^n$. Since there is one more unknown than equations, we are able to maximize the value of $u_0$, which is proportional to $\beta_2$ in every case. Summing all $n$ equations of~\eqref{eq:SCalternativesystem} leads to the nonlinear program:
\begin{equation}\label{eq:SCaltparamprob}
\begin{array}{cc}
{\rm maximize~} u_0 & \left( = \dfrac{\displaystyle\sum_{k=1}^n \left\{ \tilde{\epsilon}_k - \Im \sum_{j=1}^n u_j (x_k+{\rm i}\pi/2)^{j-1}\right\}}{\displaystyle\sum_{k=1}^n\cosh(x_k)} \right),\\
{\rm subject~to~} & h(x_k+{\rm i}\pi/2) = \tilde{\delta}_k + {\rm i}\tilde{\epsilon}_k,\quad{\rm for}\quad k=1,\ldots,n.
\end{array}
\end{equation}
Because the maximization condition is obtained by summing the constraint equations, we have one additional degree of freedom in the program~\eqref{eq:SCaltparamprob}. In order to save from premature convergence, we impose {\it ad hoc} the condition:
\begin{equation}
\left\{\begin{array}{ccc}
x_1 = 0, &~~~{\rm for} &~~~n=1,\\
|x_1+x_n|\le \bar{x}, &~~~{\rm for} &~~~n\ge2,
\end{array}\right.
\end{equation}
where $\bar{x}$ is a parameter which ensures the singularities stay reasonably close to the origin. In all our examples, we set $\bar{x} = 20$ which is sufficient. This nonlinear program is in close analogy to the Schwarz-Christoffel parameter problem. However, this method has many advantages over the Schwarz-Christoffel formula. Firstly, the map $h(t)$ is defined in terms of elementary functions and not as an integral. Secondly, the accuracy of the values of $\{u_k\}_{k=0}^n$ does not need to equal the accuracy required of the map, allowing the map~\eqref{eq:SCalternative} to be evaluated in arbitrary precision. One disadvantage of this method is that the map is an approximate solution to the original problem. Therefore, while the strip width will indeed be $2\,d = \pi$, we can expect a smaller than optimal $\beta_2$. Nevertheless, given that $\beta_2$ only has a secondary effect on the convergence rates, according to theorem~\ref{thm:DEmaxcon}, this is a small price to pay to obtain a solution method that emulates the Schwarz-Christoffel formula while being amenable to arbitrary precision calculations.

\begin{figure}[htbp]
\begin{center}
\begin{tabular}{ccc}
\includegraphics[width=0.3\textwidth]{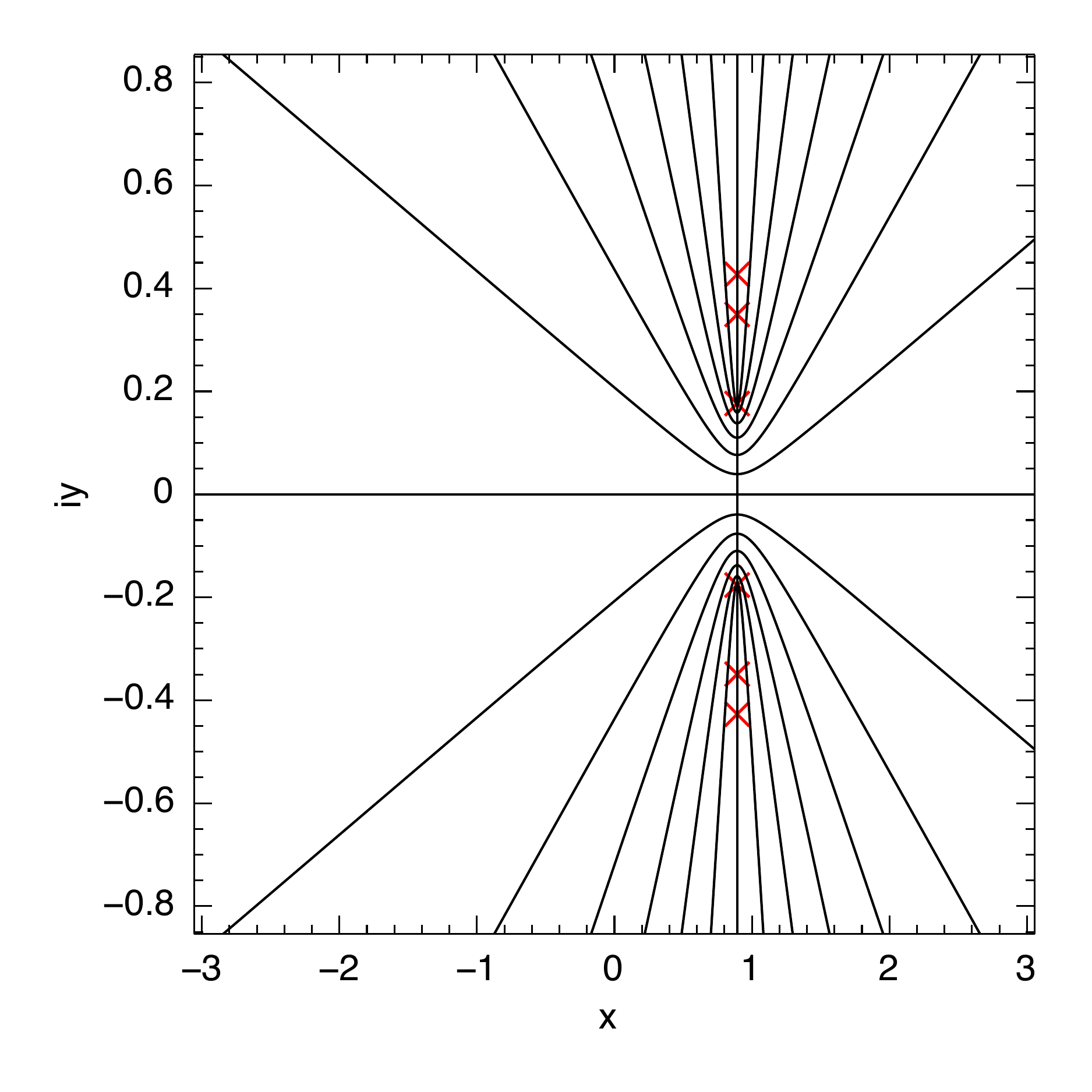}&
\includegraphics[width=0.3\textwidth]{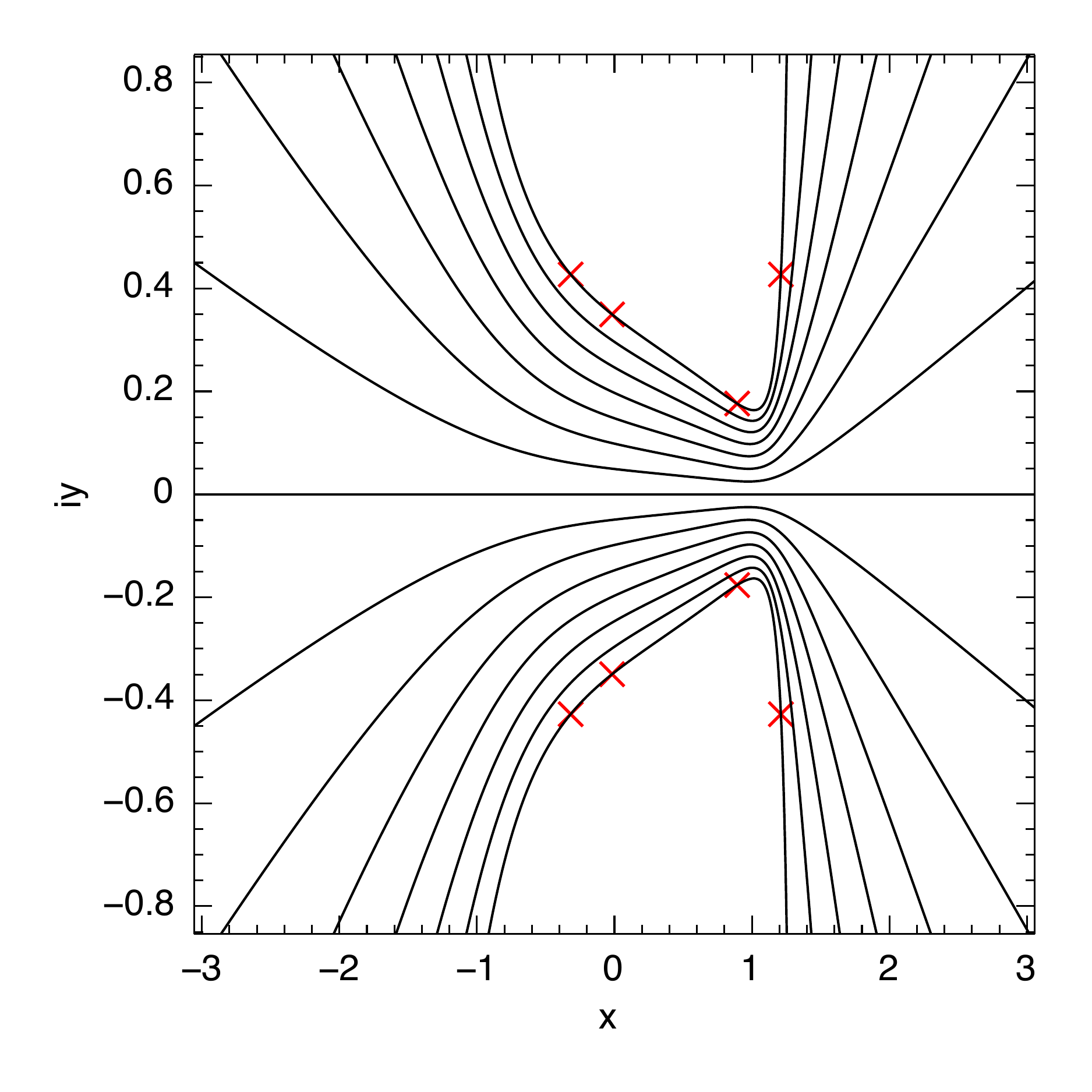}&
\includegraphics[width=0.3\textwidth]{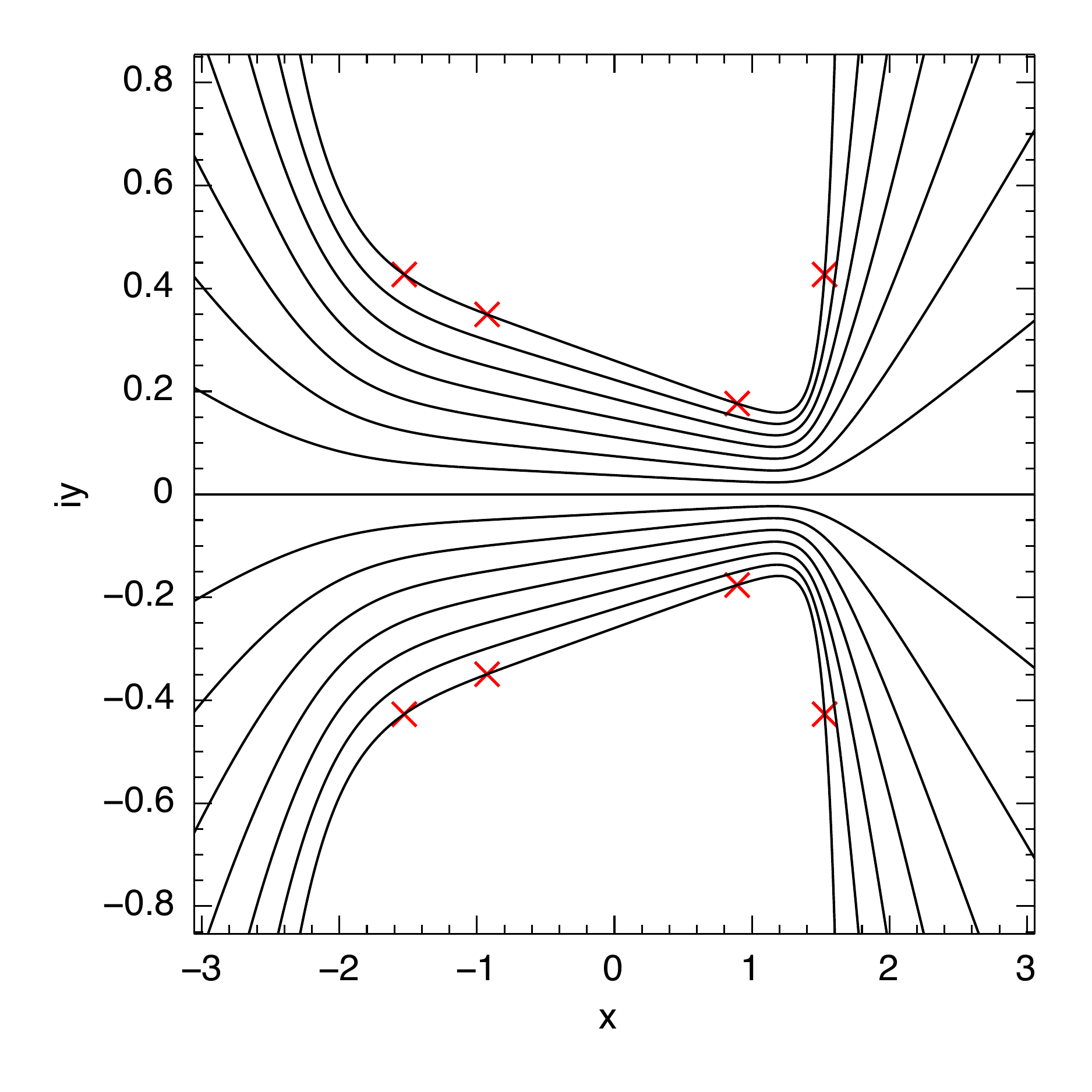}\\
(a) & (b) & (c)\\
\end{tabular}
\caption{In (a) the exact solution $\mathscr{H}(0)$, in (b) the solution $\mathscr{H}(1/2)$, and in (c) the desired solution $\mathscr{H}(1)$. An integral with these singularities is treated in Example~\ref{example:SEandDEEx2}.}
\label{fig:homotopy}
\end{center}
\end{figure}

A nonlinear program without any {\it a priori} information on the solution requires an iterative method for solving the parameter problem~\eqref{eq:SCaltparamprob}. An iterative method also requires a close initial guess to converge to the solution. To obtain an initial guess, we let $\bar{\epsilon}$ be the smallest of $\{\tilde{\epsilon}_k\}_{k=1}^n$ and $\bar{\delta}$ be the $\tilde{\delta}_k$ of the same index. Then the nonlinear program with the singularities $\{\bar{\delta}+{\rm i}\tilde{\epsilon}_k\}_{k=1}^n$ is exactly solved by:
\begin{equation}
h(t) = \bar{\epsilon}\sinh t + \bar{\delta}.
\end{equation}
A homotopy $\mathscr{H}(t)$ is then constructed between the solution with singularities $\{\bar{\delta}+{\rm i}\tilde{\epsilon}_k\}_{k=1}^n$ at $t=0$ and the solution with singularities $\{\tilde{\delta}_k+{\rm i}\tilde{\epsilon}_k\}_{k=1}^n$ at $t=1$. The interval $t\in[0,1]$ is discretized, and the nonlinear program is solved with singularities that vary linearly between the two problems and initial guesses from the solution of the previous iterate. Figure~\ref{fig:homotopy} shows this solution process.

In practice, the locations of a function's singularities may not be known in advance. This can result from either incomplete theoretical information, or a non-local nature of singularities, such as branch cuts or ``numerically singular'' terms such as the error function, which while entire, is unbounded off the real axis~\cite{Hale-Tee-31-3195-09}. In~\cite{Tee-Trefethen-28-1798-06}, an adaptive approach is taken to approximating the nearby singularities, whereby the interpolatory Chebyshev-Pad\'e approximants are constructed, and approximants' poles are taken as the loci of the singularities of the underlying function. The map is then modified to exclude these points, and the iteration of this process is the adaptive algorithm.

Because we are working with Sinc approximations, and Sinc points, we modify their algorithm to make efficient use of the information we have at hand, i.e. the Sinc sampling of the function.
\begin{definition} Let $x_k = k\,h$ be the Sinc points and let $f(x_k)$ be the $N(=2n+1)$ Sinc sampling of $f$. Then for $r+s\le 2n$, the Sinc-Pad\'e approximants $\{r/s\}_f(x)$ are given by:
\begin{equation}
\{r/s\}_f(x) = \dfrac{\displaystyle \sum_{i=0}^r p_i\,x^i}{\displaystyle 1 + \sum_{j=1}^s q_j\,x^j},
\end{equation}
where the $r+s+1$ coefficients solve the system:
\begin{equation}
\sum_{i=0}^r p_i\, x_k^i - f(x_k)\sum_{j=1}^s q_j\, x_k^j = f(x_k),
\end{equation}
for $k = -\lfloor\frac{r+s}{2}\rfloor,\ldots,\lceil\frac{r+s}{2}\rceil$.
\end{definition}

For the Chebyshev-Pad\'e approximants, the inverse cosine distribution of sample points leads to a stable linear system and the degrees of the numerator and denominator can add to equate the number of collocation points. For the Sinc-Pad\'e approximants, double exponential growth of the sample points renders the system highly ill-conditioned. Therefore, these indices must be decoupled from $n$ and the function must only be sampled near the centre. Our adaptive algorithm is based on the following principles:
\begin{enumerate}
\item Sinc-Pad\'e approximants are useful only when the Sinc approximation obtains some degree of accuracy,
\item Sinc-Pad\'e approximants are useful for $r,s = {\cal O}(\log n)$ as $n\to\infty$.
\end{enumerate}
The first principle follows from observations of our numerical experiments, and we found that a relative error of approximately $10^{-3}$ in the Sinc approximation allows for a useful Sinc-Pad\'e approximant. The second principle follows from the observation that we need not identify many singularities to remove, and that even at a logarithmic increase, the sample points tend to infinity at a single exponential rate, implying that they will ultimately cover the real line. These principles form the basis of the following algorithm.

\begin{algorithm}\label{alg:adaptiveDE}~

Set $n=1$;

{\bf while} {\rm |Relative Error|} $\ge 10^{-3}$ {\bf do}

\qquad Double $n$ and na\"ively compute the $n^{\rm th}$ double exponential approximation{\rm;}

{\bf end}{\rm;}

{\bf while} {\rm |Relative Error|} $\ge \epsilon$ {\bf do}

\qquad Compute the poles of the Sinc-Pad\'e approximant{\rm;}

\qquad Solve the nonlinear program~\eqref{eq:SCaltparamprob} for $h(t)${\rm;}

\qquad Double $n$ and compute the $n^{\rm th}$ adapted optimized approximation{\rm;}

{\bf end}{\rm.}
\end{algorithm}

\section{Examples}

In this section, we will use the proposed nonlinear program~\eqref{eq:SCaltparamprob} to maximize the convergence rate of the double exponential transformation. We compare the results of the trapezoidal rule with single, double, and optimized double exponential variable transformations on three integrals using arbitrary precision arithmetic. On a fourth integral, we use the adaptive algorithm~\ref{alg:adaptiveDE} to approximate nearby singularities.

\subsection{Example: endpoint and complex singularities}\label{example:SEandDEEx1}
We wish to evaluate the integral:
\begin{equation}\label{eq:SEandDEEx1}
\int_{-1}^1 \dfrac{\exp\left((\epsilon_1^2+(x-\delta_1)^2)^{-1}\right)\log(1-x)}{(\epsilon_2^2+(x-\delta_2)^2)\sqrt{1+x}}{\rm\,d}x = -2.04645\ldots,
\end{equation}
for the values $\delta_1 + {\rm i}\epsilon_1 = -1/2+{\rm i}$ and $\delta_2+{\rm i}\epsilon_2 = 1/2+{\rm i}/2$. This integral has two different endpoint singularities and two pairs of complex conjugate singularities of different types near the integration axis. Table~\ref{table:SEandDEEx1} summarizes the variable transformations used and the parameters in the theorems~\ref{thm:SEconvergence} and~\ref{thm:DEconvergence}.

\begin{table}[htbp]
\begin{center}
\caption{Transformations and parameters for~\eqref{eq:SEandDEEx1}.}
\label{table:SEandDEEx1}
\begin{tabular*}{\hsize}{@{\extracolsep{\fill}}c|ccc}
\hline
& Single & Double & Optimized Double \\
\hline
$\phi(t)$ & $\tanh(t/2)$ & $\tanh\left(\frac{\pi}{2}\sinh(t)\right)$ & $\tanh(h(t))$\\
$\rho$ or $\gamma$ & $1$ & $1$ & $1$\\
$\beta$ or $\beta_2$ & $1/2$ & $\pi/4$ & $0.06956$\\
$d$ & $1.10715$ & $0.34695$ & $\pi/2$\\
\hline
\end{tabular*}
\end{center}
\end{table}

\begin{figure}[htbp]
\begin{center}
\begin{tabular}{ccc}
\includegraphics[width=0.3\textwidth]{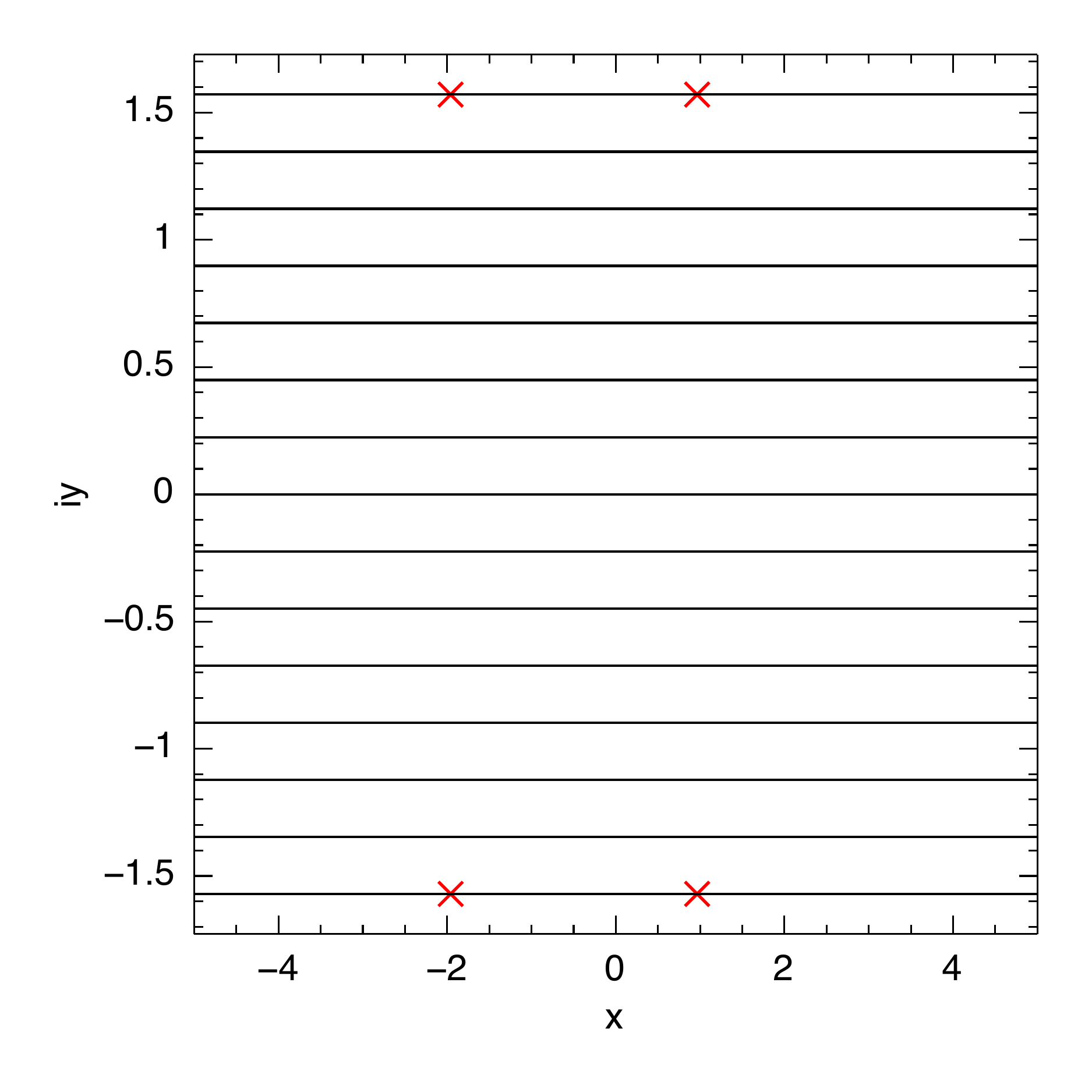}&
\includegraphics[width=0.3\textwidth]{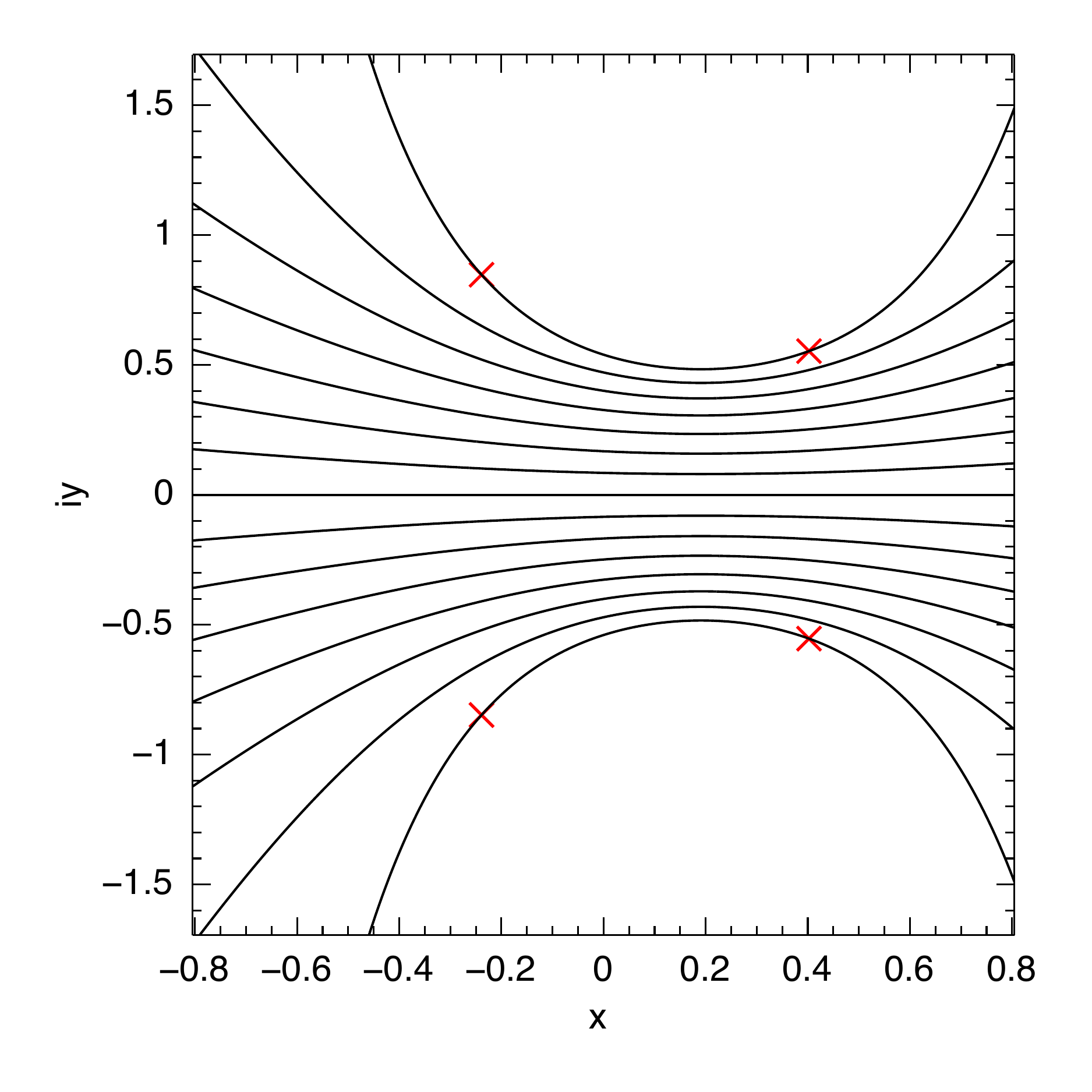}&
\includegraphics[width=0.3\textwidth]{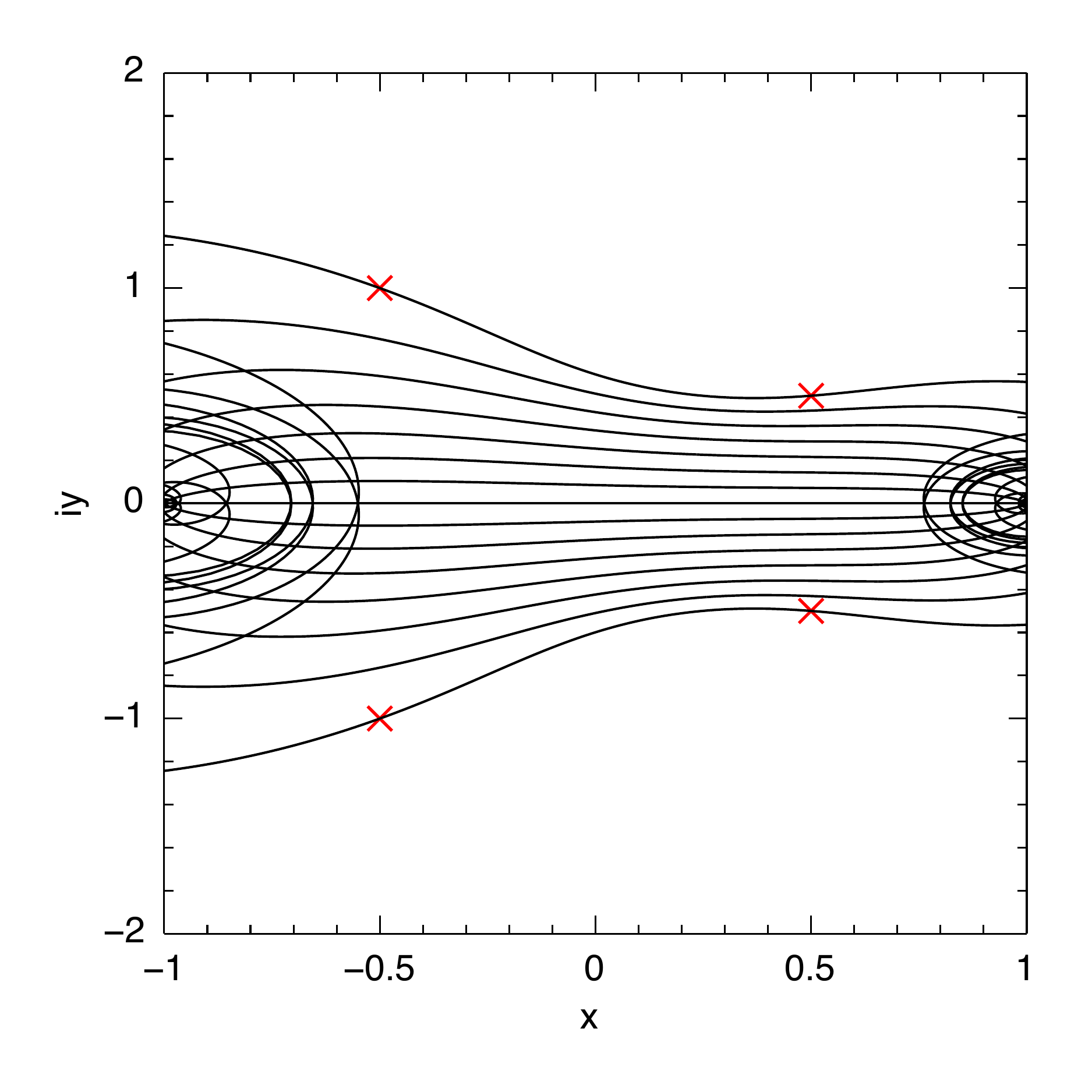}\\
(a) & (b) & (c)\\
\end{tabular}
\caption{In (a) the plot of the strip $\mathscr{D}_{\frac{\pi}{2}}$ with singularities located on the boundary, in (b) the optimized map $h(\cdot)$, and in (c) the optimized DE map. In all three cases, the crosses track the singularities.}
\label{fig:Example1maps}
\end{center}
\end{figure}

In addition, the optimized transformation is given by:
\begin{equation}
h(t) \approx 0.13912\sinh(t) + 0.19081 + 0.21938\,t.
\end{equation}
Figure~\ref{fig:Example1maps} shows the three stages of the optimized double exponential map.

\begin{figure}[htbp]
\begin{center}
\begin{tabular}{cc}
\includegraphics[width=0.45\textwidth]{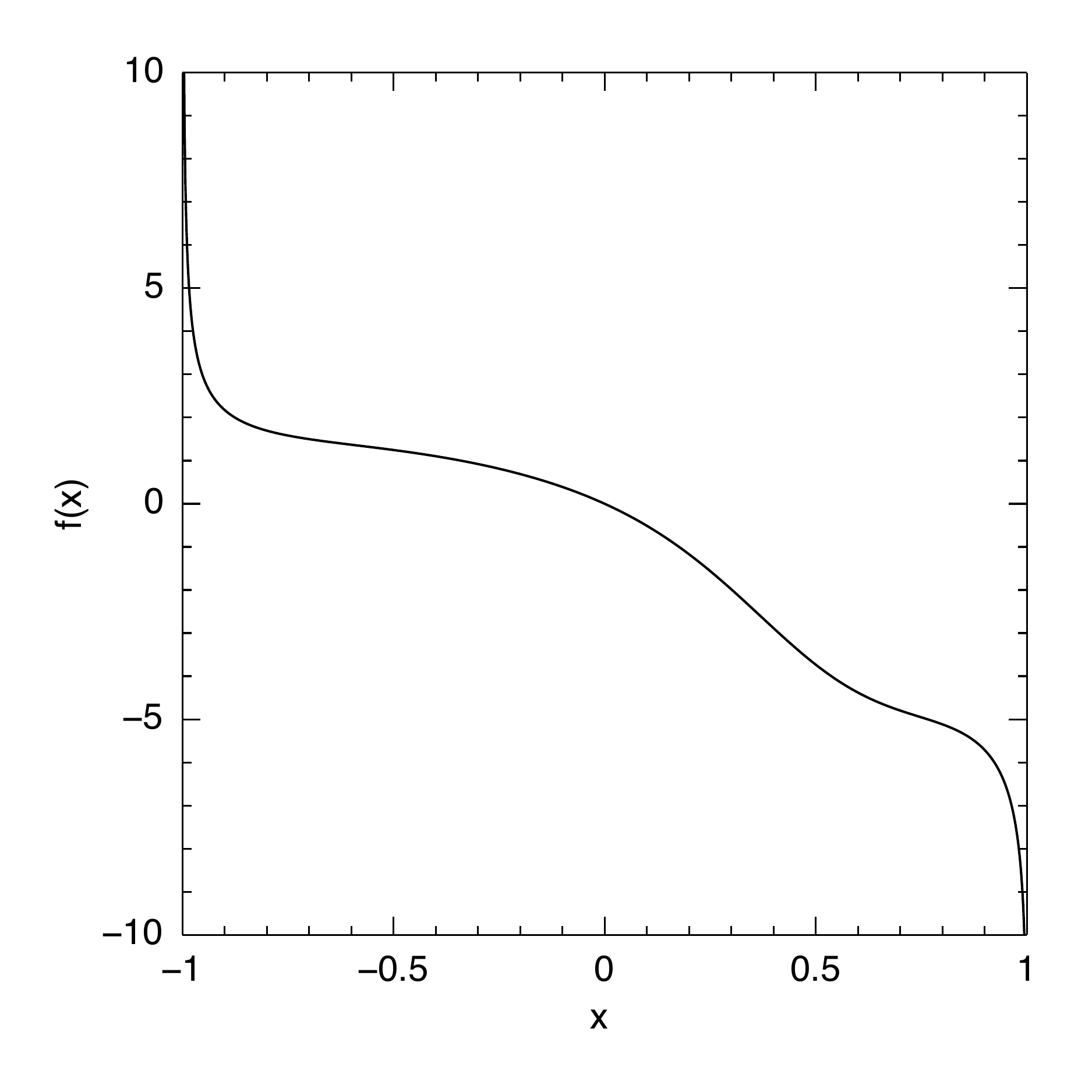}&
\includegraphics[width=0.45\textwidth]{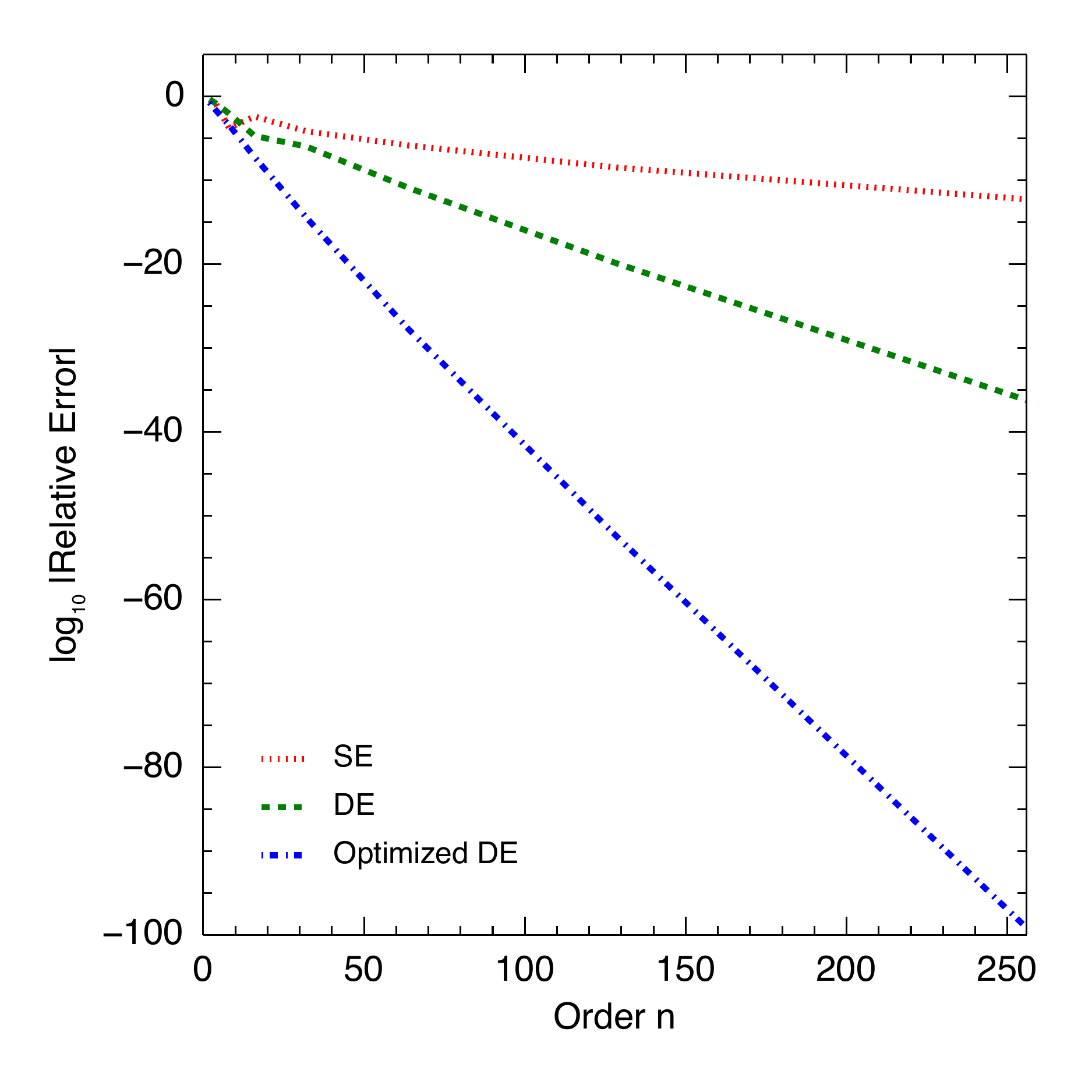}\\
(a) & (b)\\
\end{tabular}
\caption{In (a) the plot of the integrand of~\eqref{eq:SEandDEEx1} and in (b) the performance of the trapezoidal rule with single, double, and optimized double exponential variable transformations.}
\label{fig:Example1ploterr}
\end{center}
\end{figure}

In Figure~\ref{fig:Example1ploterr} (a), the integrand of~\eqref{eq:SEandDEEx1} is shown, and in Figure~\ref{fig:Example1ploterr} (b), the logarithm of the relative errors of the trapezoidal rule of order $n$ with single, double, and optimized double exponential variable transformations are plotted. The increase in convergence rate using the optimized variable transformation is a significant increase in efficiency over the double exponential transformation.

\subsection{Example: eight different complex conjugate singularities}\label{example:SEandDEEx2} We wish to evaluate the integral:
\begin{equation}\label{eq:SEandDEEx2}
\int_{-\infty}^{+\infty} \dfrac{\exp\left(10(\epsilon_1^2+(x-\delta_1)^2)^{-1}\right)\cos\left(10(\epsilon_2^2+(x-\delta_2)^2)^{-1}\right)}{(\epsilon_3^2+(x-\delta_3)^2)\sqrt{\epsilon_4^2+(x-\delta_4)^2}}{\rm d}x = 15.01336\ldots,
\end{equation}
for the values $\delta_1 + {\rm i}\epsilon_1 = -2+{\rm i}$, $\delta_2+{\rm i}\epsilon_2 = -1+{\rm i}/2$, $\delta_3+{\rm i}\epsilon_3 = 1+{\rm i}/4$, and $\delta_4+{\rm i}\epsilon_4 = 2 +{\rm i}$.
Table~\ref{table:SEandDEEx2} summarizes the variable transformations used and the parameters in the theorems~\ref{thm:SEconvergence} and~\ref{thm:DEconvergence}.

\begin{table}[htbp]
\begin{center}
\caption{Transformations and parameters for~\eqref{eq:SEandDEEx2}.}
\label{table:SEandDEEx2}
\begin{tabular*}{\hsize}{@{\extracolsep{\fill}}c|ccc}
\hline
& Single & Double & Optimized Double \\
\hline
$\phi(t)$ & $\sinh(t)$ & $\sinh\left(\frac{\pi}{2}\sinh(t)\right)$ & $\sinh(h(t))$\\
$\rho$ or $\gamma$ & $1$ & $1$ & $1$\\
$\beta$ or $\beta_2$ & $2$ & $\pi/2$ & $5.7715\times10^{-6}$\\
$d$ & $0.35260$ & $0.22640$ & $\pi/2$\\
\hline
\end{tabular*}
\end{center}
\end{table}

\begin{figure}[htbp]
\begin{center}
\begin{tabular}{ccc}
\includegraphics[width=0.3\textwidth]{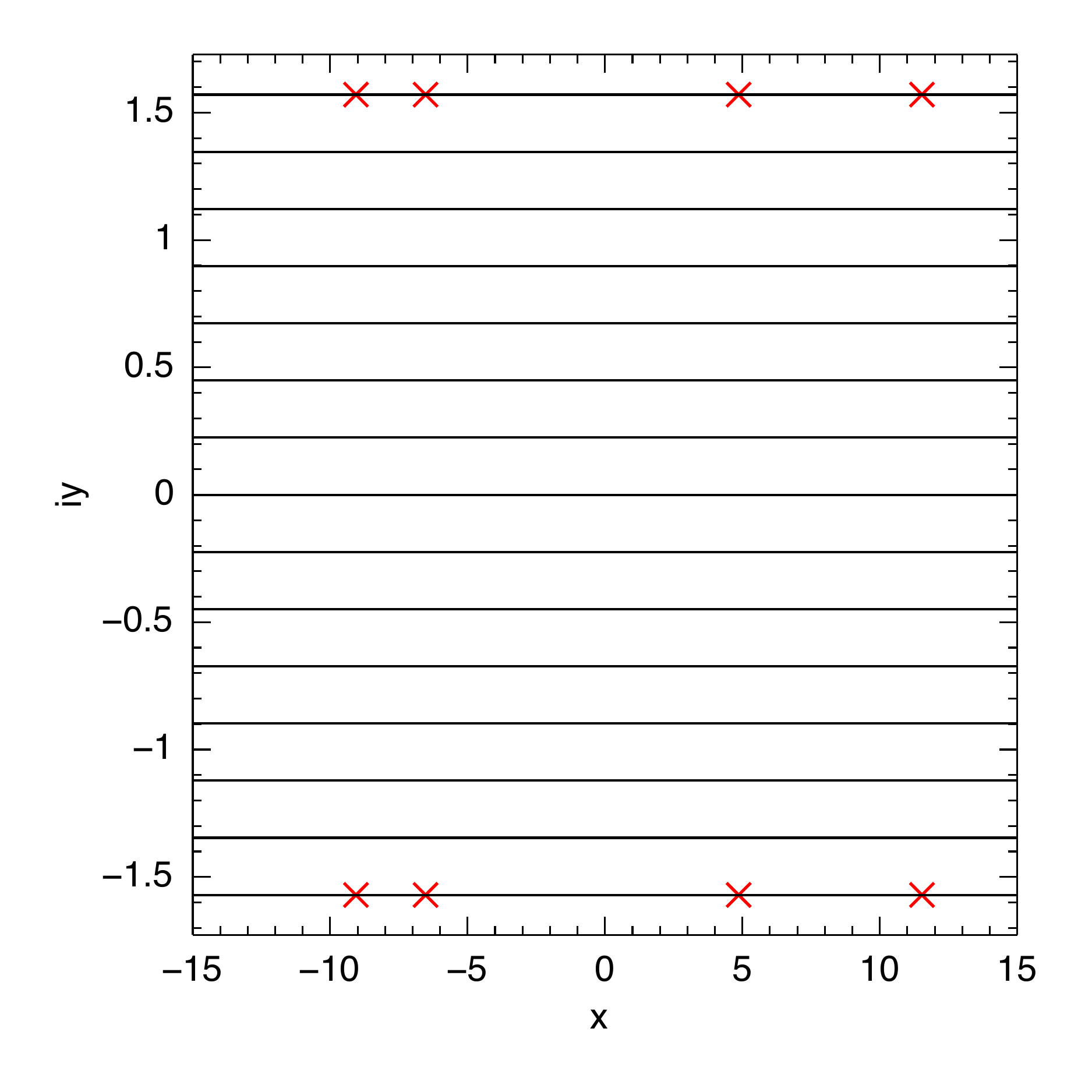}&
\includegraphics[width=0.3\textwidth]{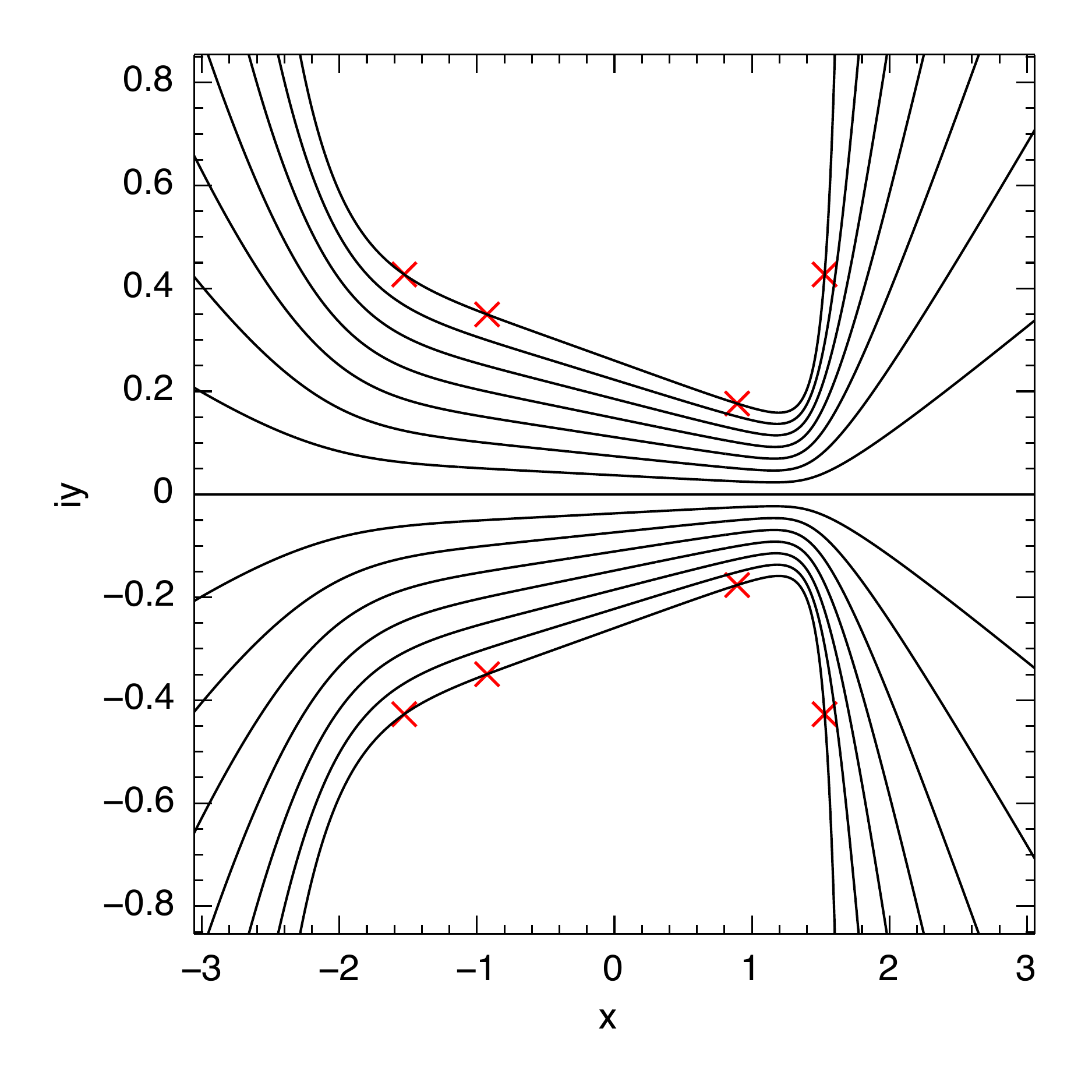}&
\includegraphics[width=0.3\textwidth]{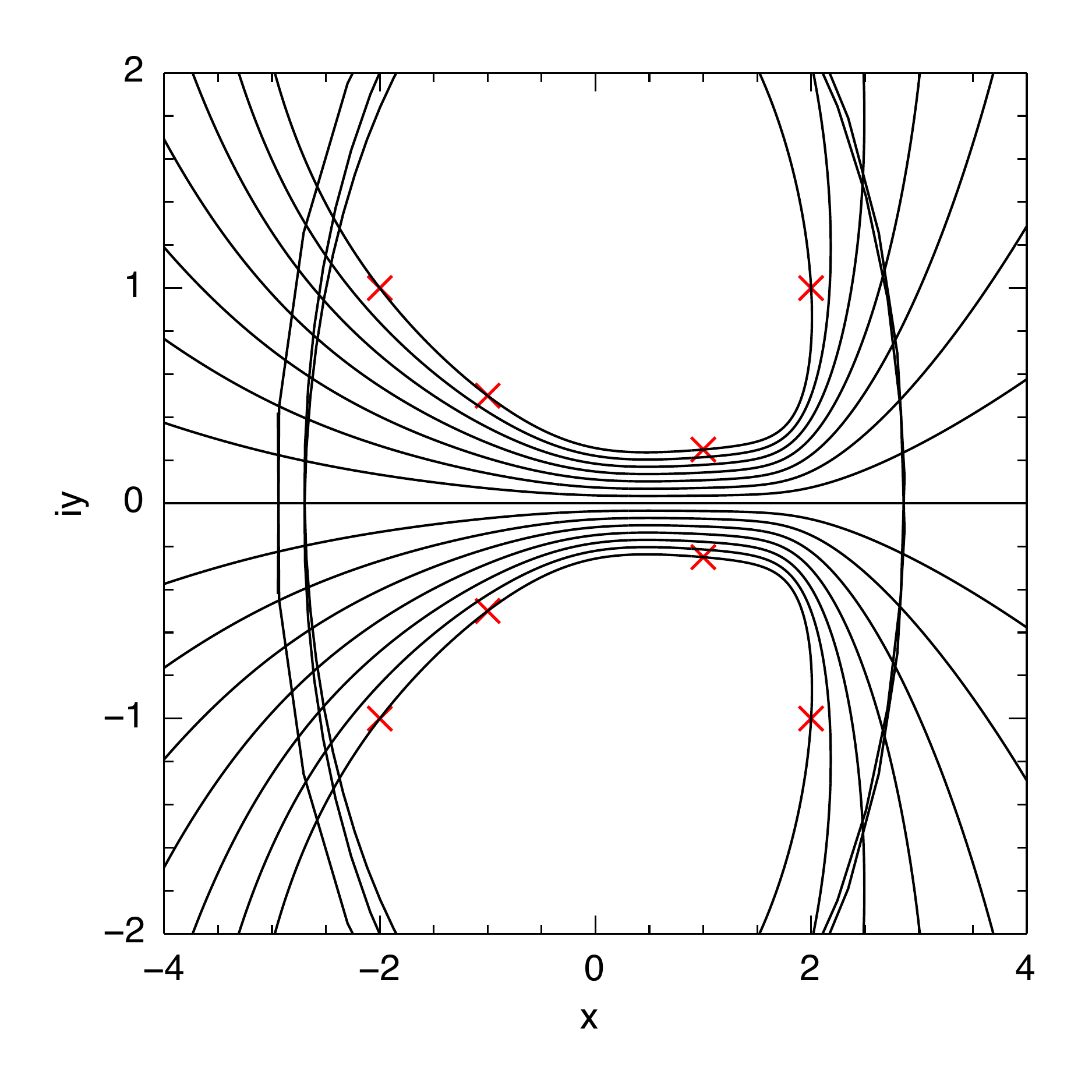}\\
(a) & (b) & (c)\\
\end{tabular}
\caption{In (a) the plot of the strip $\mathscr{D}_{\frac{\pi}{2}}$ with singularities located on the boundary, in (b) the optimized map $h(\cdot)$, and in (c) the optimized DE map. In all three cases, the crosses track the singularities.}
\label{fig:Example2maps}
\end{center}
\end{figure}

In addition, the optimized transformation is given by:
\begin{align}
h(t) \approx &~ 5.7715\times10^{-6}\sinh(t) + 0.25431 + 0.14936\,t\nonumber\\
& - 4.5433\times10^{-3}\,t^2 + 9.9880\times10^{-5}\,t^3.
\end{align}
Figure~\ref{fig:Example2maps} shows the three stages of the optimized double exponential map.

\begin{figure}[htbp]
\begin{center}
\begin{tabular}{cc}
\includegraphics[width=0.45\textwidth]{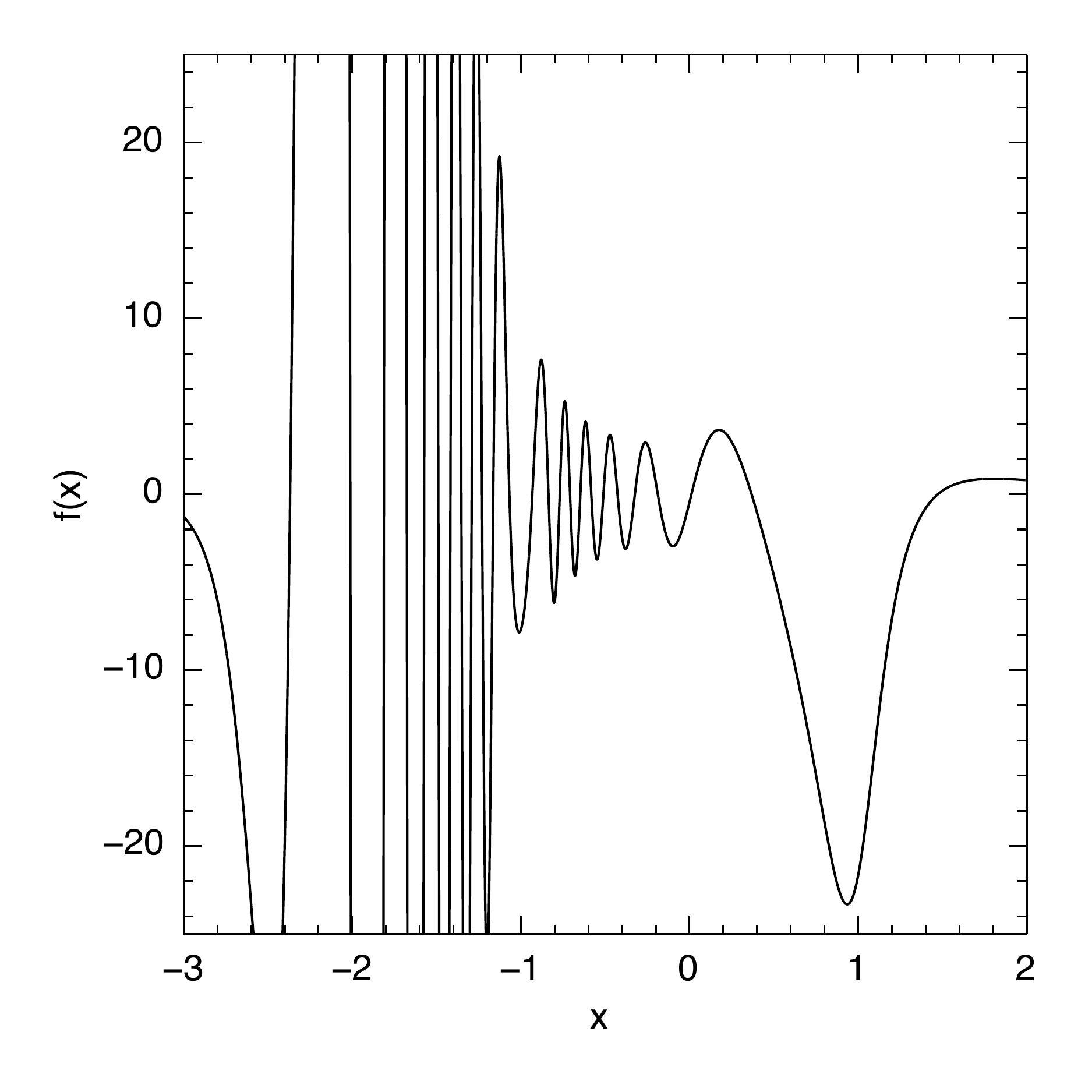}&
\includegraphics[width=0.45\textwidth]{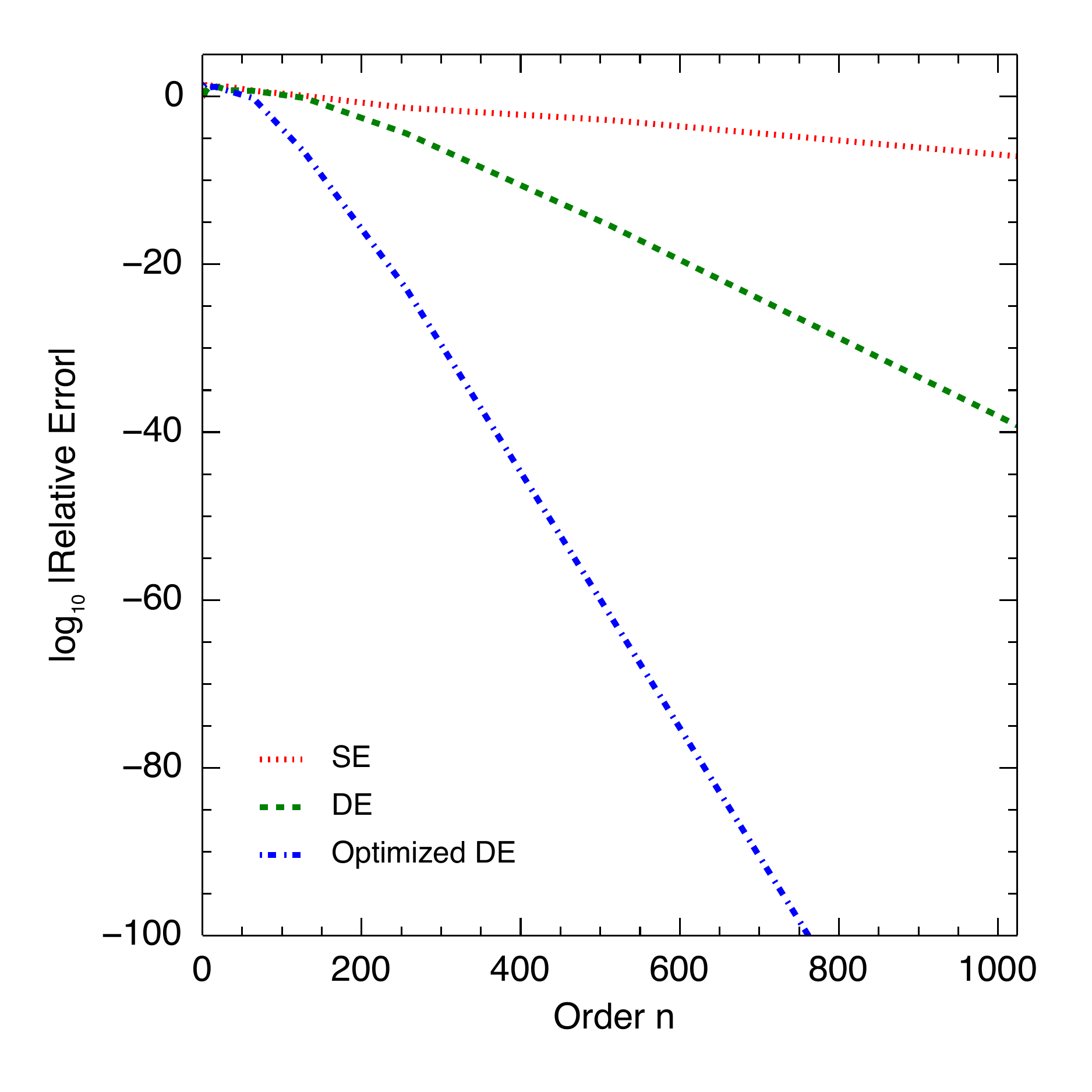}\\
(a) & (b)\\
\end{tabular}
\caption{In (a) the plot of the integrand of~\eqref{eq:SEandDEEx2} and in (b) the performance of the trapezoidal rule with single, double, and optimized double exponential variable transformations.}
\label{fig:Example2ploterr}
\end{center}
\end{figure}

In Figure~\ref{fig:Example2ploterr} (a), the integrand of~\eqref{eq:SEandDEEx2} is shown, and in Figure~\ref{fig:Example2ploterr} (b), the logarithm of the relative errors of the trapezoidal rule of order $n$ with single, double, and optimized double exponential variable transformations are plotted. The increase in convergence rate using the optimized variable transformation is a significant increase in efficiency over the double exponential transformation.

\subsection{Example: for Goursat's infinite integral}\label{example:SEandDEEx3} We wish to evaluate the integral:
\begin{equation}\label{eq:SEandDEEx3}
\int_{0}^{+\infty} \dfrac{x{\rm\,d}x}{1+x^6\sinh^2x} = 0.50368\ldots,
\end{equation}
which is evaluated in~\cite{Hatano-Ninomiya-Sugiura-Hasegawa-52-213-09,Ooura-249-1-13} as part of a high precision numerical evaluation of Goursat's infinite integral. While there are an infinite number of poles in the complex plane due to the $\sinh$ function, a four-parameter solution $h(t)$ can be found using the nearest poles, while excluding the remainder. This shows the incredible versatility of the proposed optimization approach, because the same optimal asymptotic convergence rate is obtained in the complicated situation of an infinite number of singularities while not leading to a more complicated solution process. Table~\ref{table:SEandDEEx3} summarizes the variable transformations used and the parameters in the theorems~\ref{thm:SEconvergence} and~\ref{thm:DEconvergence}. For the sake of comparison, we use the same double exponential transformation used in~\cite{Hatano-Ninomiya-Sugiura-Hasegawa-52-213-09}.

\begin{table}[htbp]
\begin{center}
\caption{Transformations and parameters for~\eqref{eq:SEandDEEx3}.}
\label{table:SEandDEEx3}
\begin{tabular*}{\hsize}{@{\extracolsep{\fill}}c|ccc}
\hline
& Single & Double & Optimized Double \\
\hline
$\phi(t)$ & $\log(e^t+1)$ & $\exp(0.22t-0.017e^{-t})$ & $\log(e^{h(t)}+1)$\\
$\rho$ or $\gamma$ & $1$ & $0.22$ & $1$\\
$\beta$ or $\beta_2$ & $2$ & $2$ & $0.26725$\\
$d$ & $1.13615$ & $1.58223$ & $\pi/2$\\
\hline
\end{tabular*}
\end{center}
\end{table}

\begin{figure}[htbp]
\begin{center}
\begin{tabular}{ccc}
\includegraphics[width=0.3\textwidth]{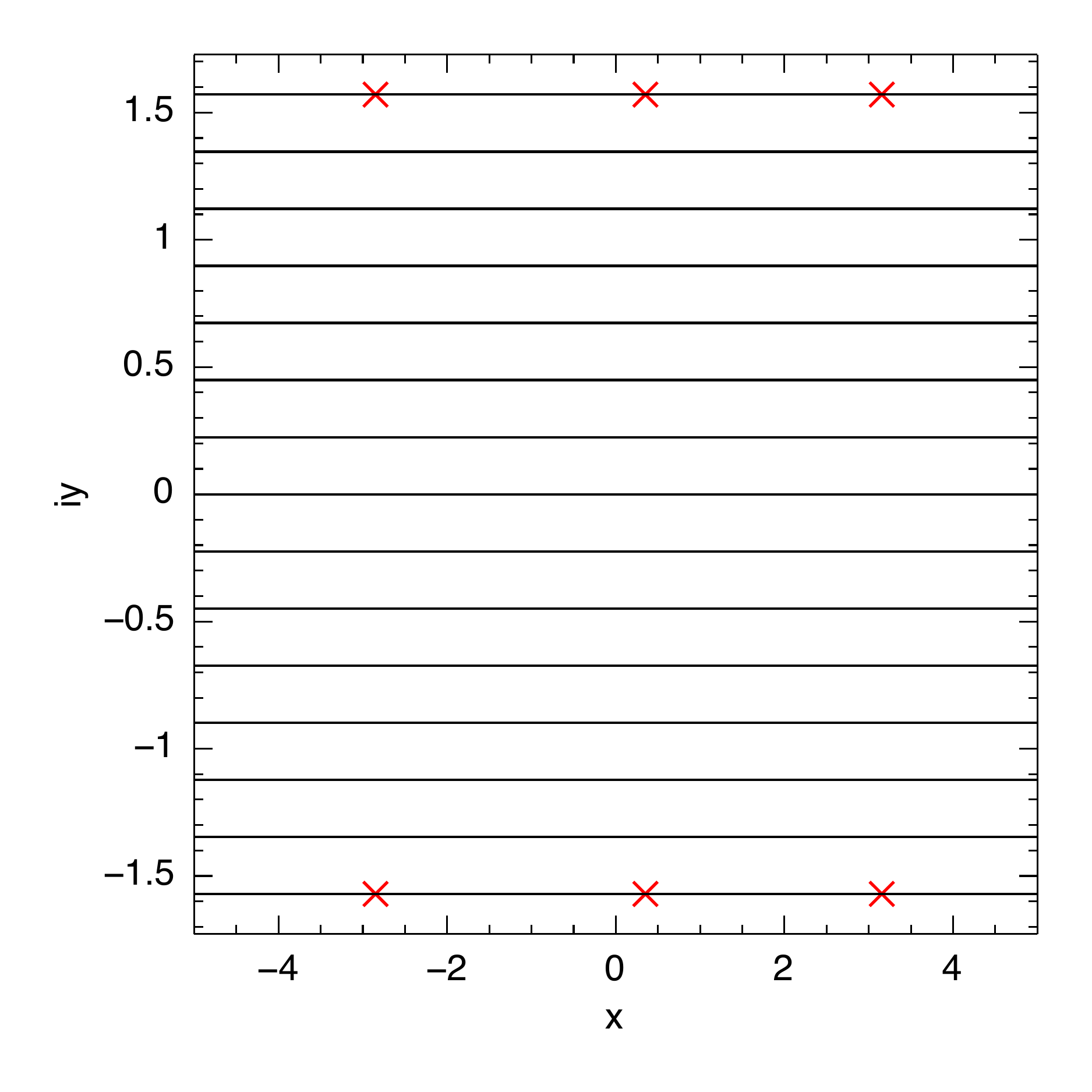}&
\includegraphics[width=0.3\textwidth]{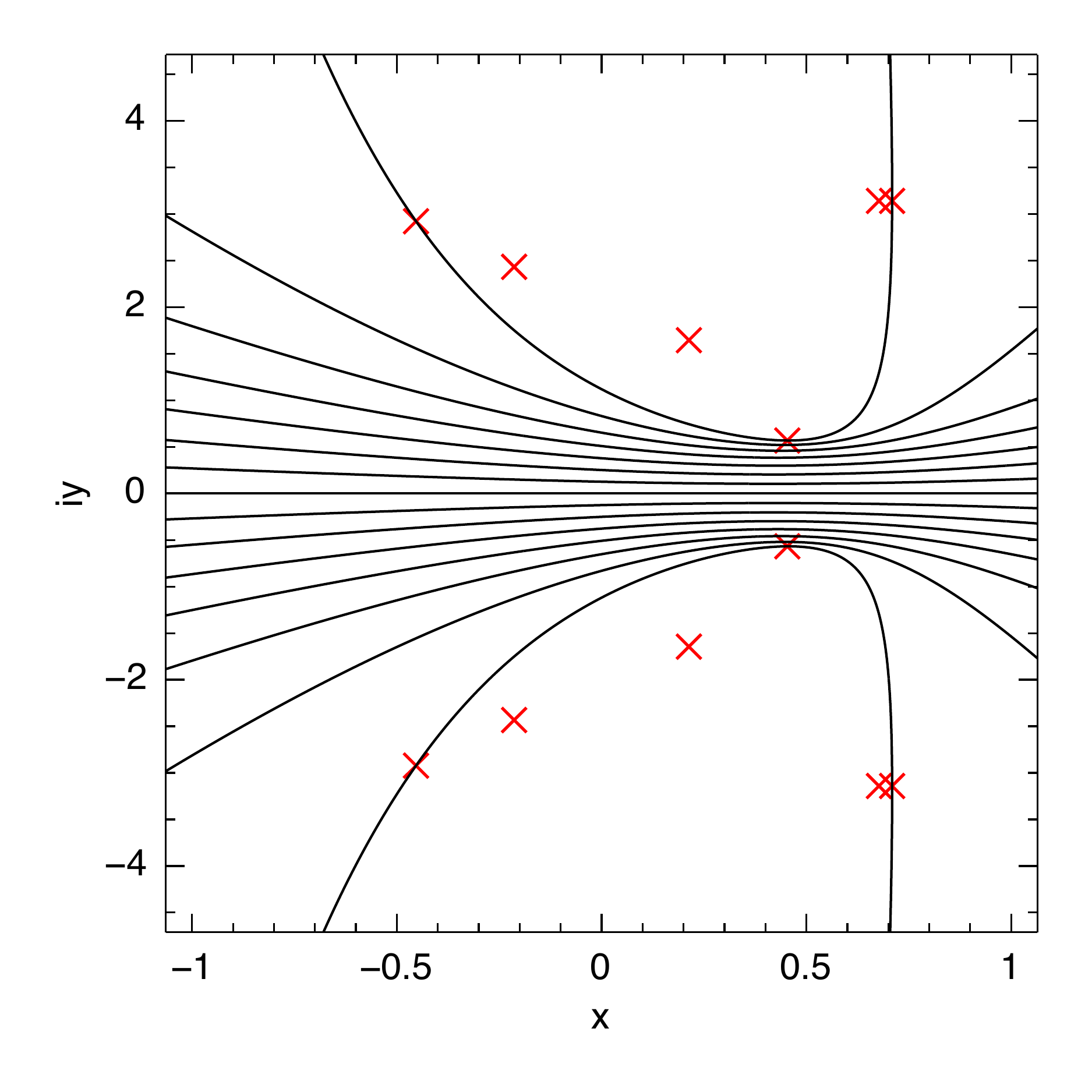}&
\includegraphics[width=0.3\textwidth]{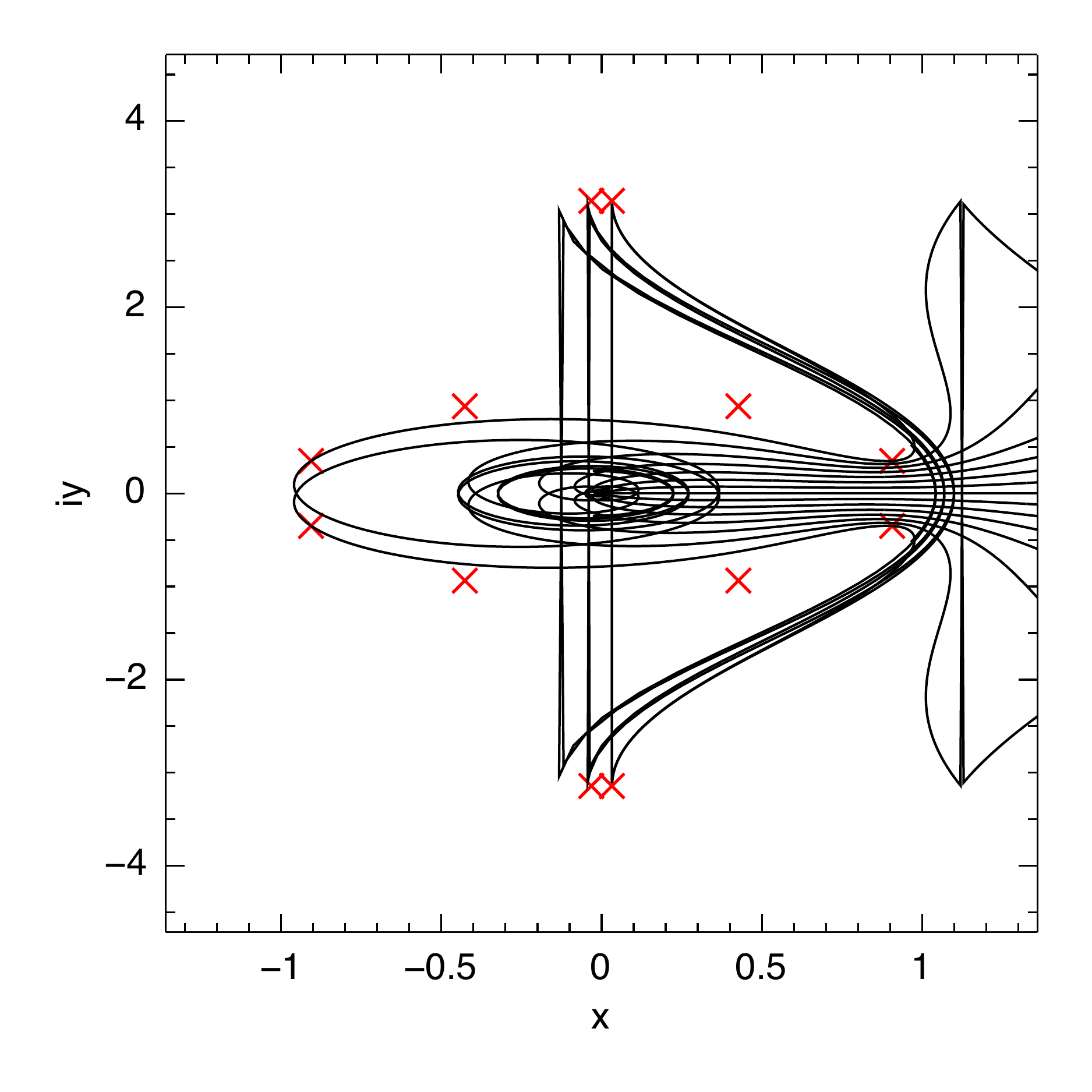}\\
(a) & (b) & (c)\\
\end{tabular}
\caption{In (a) the plot of the strip $\mathscr{D}_{\frac{\pi}{2}}$ with singularities located on the boundary, in (b) the optimized map $h(\cdot)$, and in (c) the optimized DE map. In all three cases, the crosses track the singularities.}
\label{fig:Example3maps}
\end{center}
\end{figure}

In addition, the optimized transformation is given by:
\begin{equation}
h(t) \approx 0.26725\sinh(t) + 0.30707 + 0.20337\,t - 0.031966\,t^2.
\end{equation}
Figure~\ref{fig:Example3maps} shows the three stages of the optimized double exponential map.

\begin{figure}[htbp]
\begin{center}
\begin{tabular}{cc}
\includegraphics[width=0.45\textwidth]{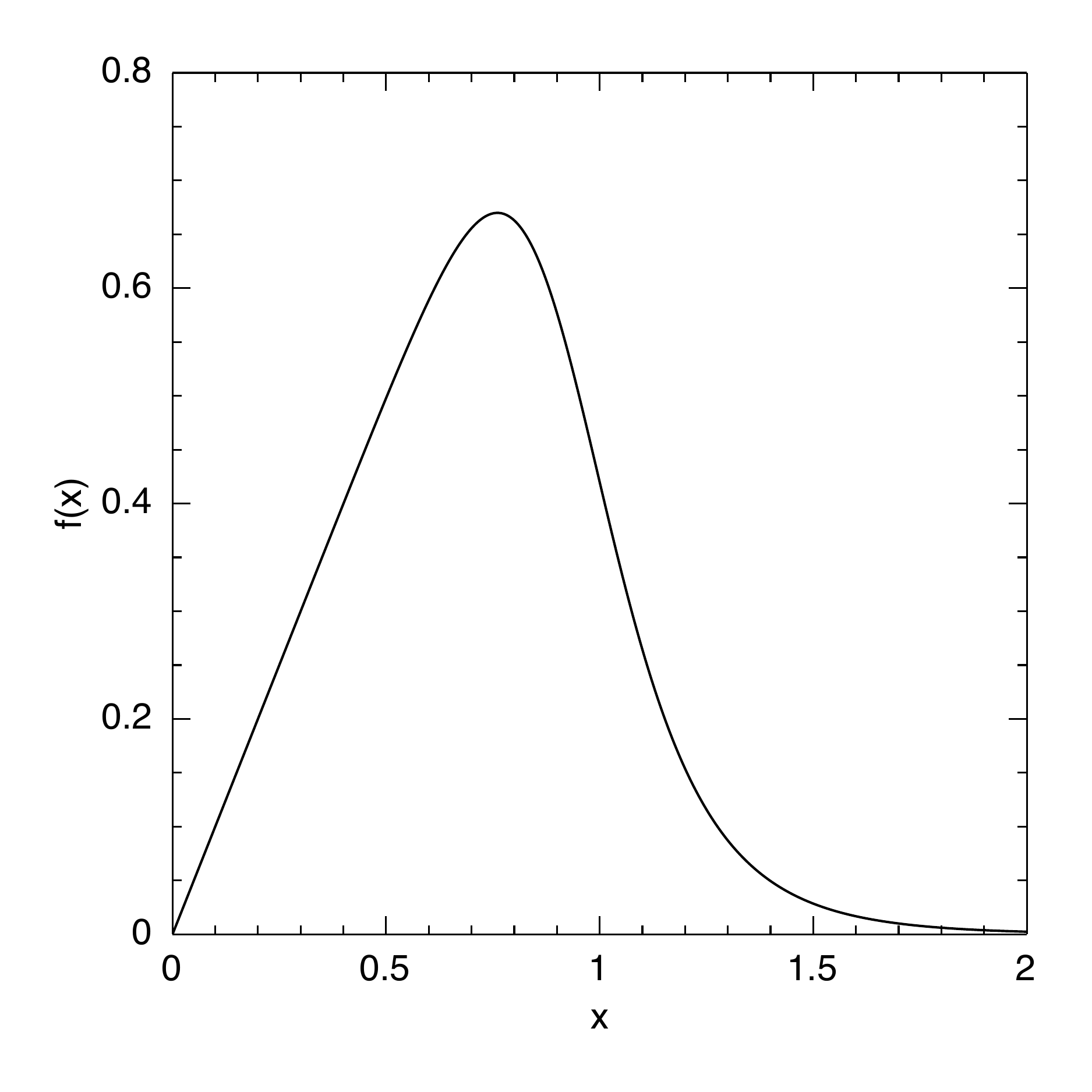}&
\includegraphics[width=0.45\textwidth]{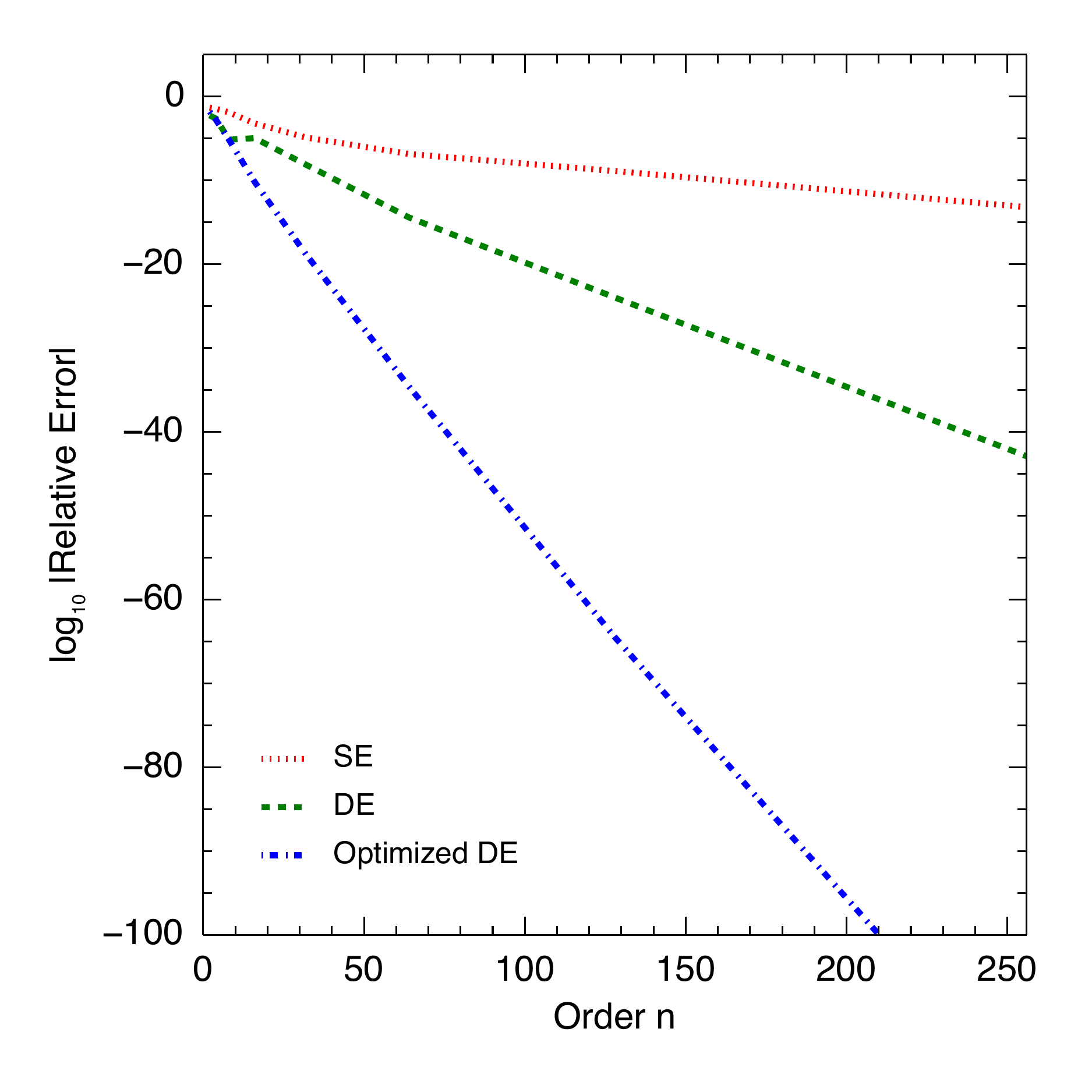}\\
(a) & (b)\\
\end{tabular}
\caption{In (a) the plot of the integrand of~\eqref{eq:SEandDEEx3} and in (b) the performance of the trapezoidal rule with single, double, and optimized double exponential variable transformations.}
\label{fig:Example3ploterr}
\end{center}
\end{figure}

In Figure~\ref{fig:Example3ploterr} (a), the integrand of~\eqref{eq:SEandDEEx3} is shown, and in Figure~\ref{fig:Example3ploterr} (b), the logarithm of the relative errors of the trapezoidal rule of order $n$ with single, double, and optimized double exponential variable transformations are plotted. The increase in convergence rate using the optimized variable transformation is a significant increase in efficiency over the double exponential transformation.

In~\cite{Hatano-Ninomiya-Sugiura-Hasegawa-52-213-09}, the authors obtained the relative error of $10^{-72}$ with $480$ function evaluations, and claim this to be a nearly minimal number. However, as can be seen in Figure~\ref{fig:Example3ploterr} (b), the optimized double exponential transformation can obtain the same relative error with only approximately $140$ function evaluations. Neither of these results, however, compares to the million-digit algorithm of~\cite{Ooura-249-1-13} using the Hilbert transform to construct a conjugate function, thereby removing all the singularities.

\subsection{Example: adaptive optimization via Sinc-Pad\'e approximants}\label{example:SEandDEEx4} In this example, we will show the performance of the adaptive optimized method using the Sinc-Pad\'e approximants to approximately locate the singularities. We wish to evaluate the integral:
\begin{equation}\label{eq:SEandDEEx4}
\int_0^{\infty} \dfrac{x{\rm\,d}x}{\sqrt{\epsilon_1^2+(x-\delta_1)^2}(\epsilon_2^2+(x-\delta_2)^2)(\epsilon_3^2+(x-\delta_3)^2)} = 12.55613\ldots,
\end{equation}
for the values $\delta_1 + {\rm i}\epsilon_1 = 1+{\rm i}$, $\delta_2+{\rm i}\epsilon_2 = 2+{\rm i}/2$, and $\delta_3+{\rm i}\epsilon_3 = 3+{\rm i}/3$. Table~\ref{table:SEandDEEx4} summarizes the variable transformations used and the parameters in the theorems~\ref{thm:SEconvergence} and~\ref{thm:DEconvergence}.

\begin{table}[htbp]
\begin{center}
\caption{Transformations and parameters for~\eqref{eq:SEandDEEx4}.}
\label{table:SEandDEEx4}
\begin{tabular*}{\hsize}{@{\extracolsep{\fill}}c|ccc}
\hline
& Single & Double & Optimized Double \\
\hline
$\phi(t)$ & $\exp(t)$ & $\exp\left(\frac{\pi}{2}\sinh(t)\right)$ & $\exp(h(t))$\\
$\rho$ or $\gamma$ & $1$ & $1$ & $1$\\
$\beta$ or $\beta_2$ & $2$ & $\pi/2$ & $9.4353\times10^{-3}$\\
$d$ & $0.11066$ & $0.05762$ & $\pi/2$\\
\hline
\end{tabular*}
\end{center}
\end{table}

Table~\ref{table:SEandDEEx4SincPade} shows the evolution of the six nearest roots of the Sinc-Pad\'e approximants to the integration contour. The degrees of the Sinc-Pad\'e approximants increase as $r = \log_2(n) -2$ and $s = \log_2(n) +2$.

\begin{table}[htbp]
\begin{center}
\caption{Evolution of the six nearest roots of the Sinc-Pad\'e approximants.}
\label{table:SEandDEEx4SincPade}
\begin{tabular*}{\hsize}{@{\extracolsep{\fill}}c|ccc}
\sphline
$n$ & $\delta_1\pm{\rm i}\epsilon_1$ & $\delta_2\pm{\rm i}\epsilon_2$ & $\delta_3\pm{\rm i}\epsilon_3$\\
\sphline
$2^5$ & $-0.58097 \pm 1.2106{\rm i}$ & $2.0717 \pm 0.28089{\rm i}$ & $3.0298 \pm 0.45170{\rm i}$\\
$2^6$ & $-0.18822 \pm 1.3571{\rm i}$ & $2.0008 \pm 0.49777{\rm i}$ & $3.0004 \pm 0.33311{\rm i}$\\
$2^7$ & ~~\,$0.14091 \pm 1.3982{\rm i}$ & $1.9963 \pm 0.48734{\rm i}$ & $3.0009 \pm 0.33279{\rm i}$\\
$2^8$ & ~~\,$0.41762 \pm 1.3767{\rm i}$ & $2.0498 \pm 0.39481{\rm i}$ & $3.0022 \pm 0.35279{\rm i}$\\
\sphline
Exact & ~~\,$1.0000 \pm 1.0000{\rm i}$ & $2.0000\pm0.50000{\rm i}$ & $3.0000\pm0.3333\bar{3}{\rm i}$\\
\hline
\end{tabular*}
\end{center}
\end{table}

Table~\ref{table:SEandDEEx4adapt} shows the evolution of the adaptive map. The coefficients of the optimized map are also shown for comparison.

\begin{table}[htbp]
\begin{center}
\caption{Evolution of the coefficients of the adaptive map.}
\label{table:SEandDEEx4adapt}
\begin{tabular*}{\hsize}{@{\extracolsep{\fill}}c|cccc}
\sphline
$n$ & $u_0$ & $u_1$ & $u_2$ & $u_3$\\
\sphline
$2^5$ & $3.1344\times10^{-3}$ & $0.88233$ & $0.072018$ & $-1.6222\times10^{-3}$\\
$2^6$ & $9.1841\times10^{-3}$ & $0.95544$ & $0.073207$ & $-7.1021\times10^{-3}$\\
$2^7$ & $8.6135\times10^{-3}$ & $0.95359$ & $0.072918$ & $-6.7917\times10^{-3}$\\
$2^8$ & $6.1605\times10^{-3}$ & $0.93730$ & $0.071916$ & $-4.7927\times10^{-3}$\\
\sphline
Optimized & $9.4353\times10^{-3}$ & $0.93351$ & $0.084087$ & $-9.9846\times10^{-3}$\\
\hline
\end{tabular*}
\end{center}
\end{table}

\begin{figure}[htbp]
\begin{center}
\begin{tabular}{ccc}
\includegraphics[width=0.3\textwidth]{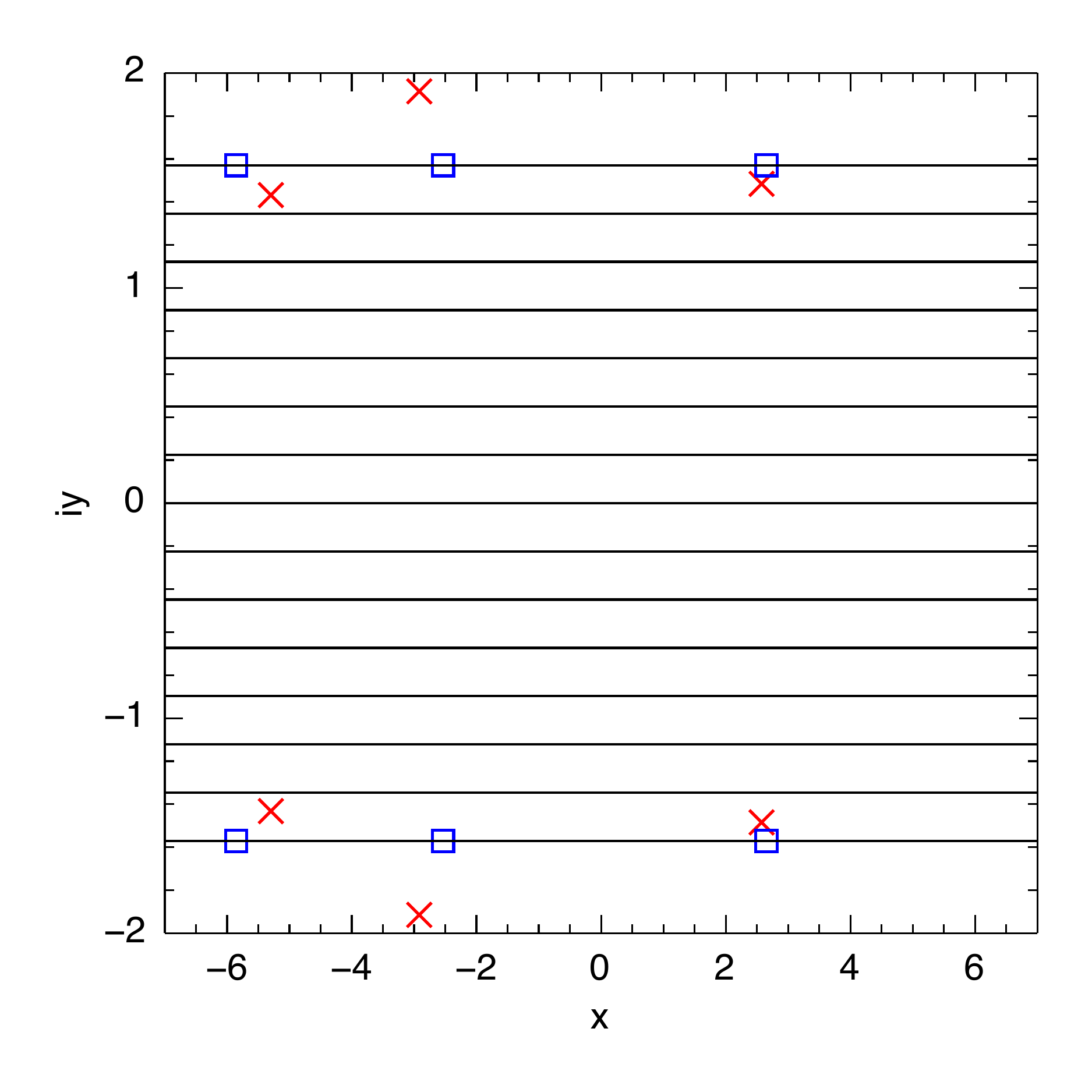}&
\includegraphics[width=0.3\textwidth]{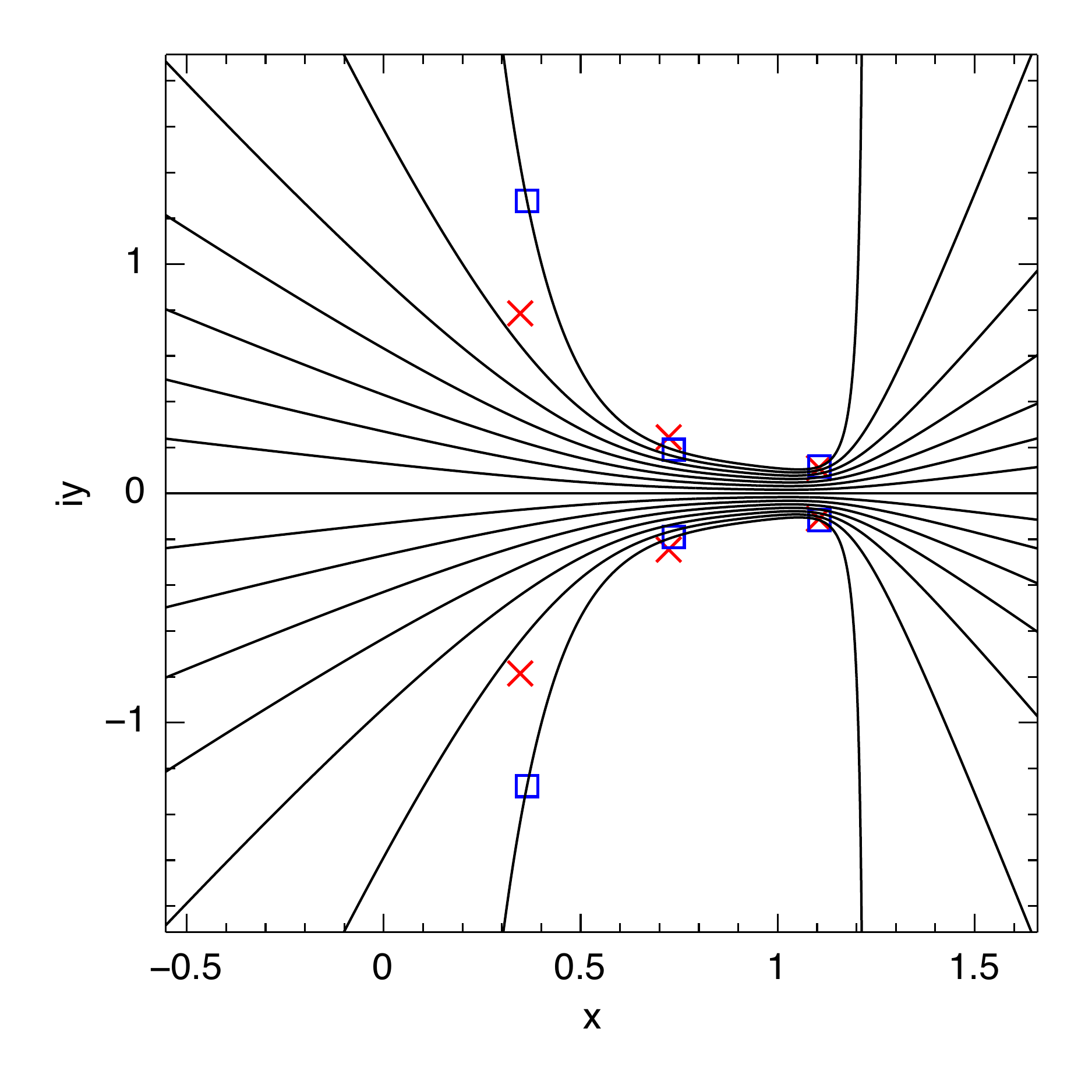}&
\includegraphics[width=0.3\textwidth]{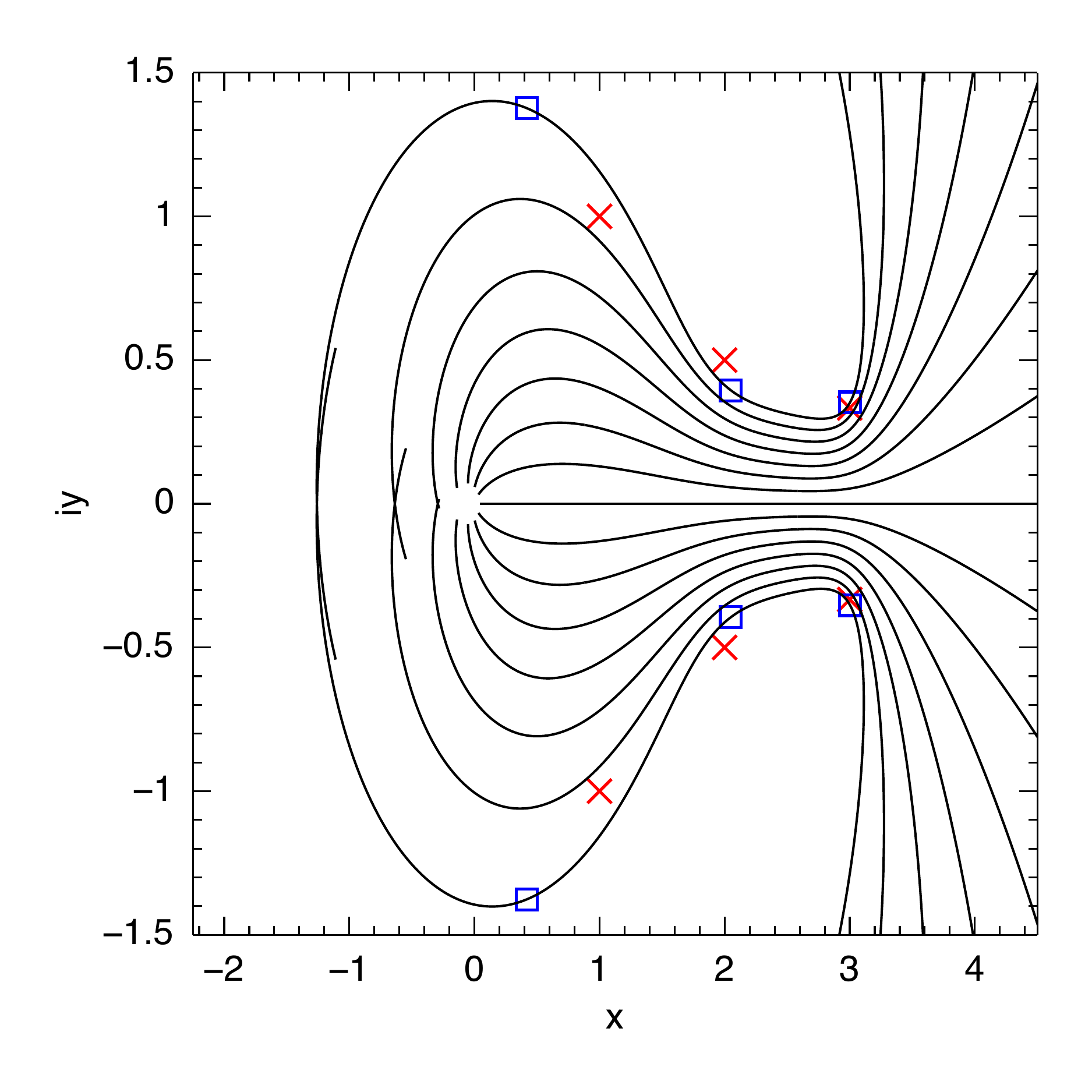}\\
(a) & (b) & (c)\\
\end{tabular}
\caption{In (a) the plot of the strip $\mathscr{D}_{\frac{\pi}{2}}$ with singularities located on the boundary, in (b) the adaptive map $h(\cdot)$, and in (c) the optimized DE map. In all three cases, the crosses track the singularities and the squares track the roots of the Sinc-Pad\'e approximant for $n=2^8$.}
\label{fig:Example4maps}
\end{center}
\end{figure}

Figure~\ref{fig:Example4maps} shows the three stages of the adaptive double exponential map.

\begin{figure}[htbp]
\begin{center}
\begin{tabular}{cc}
\includegraphics[width=0.45\textwidth]{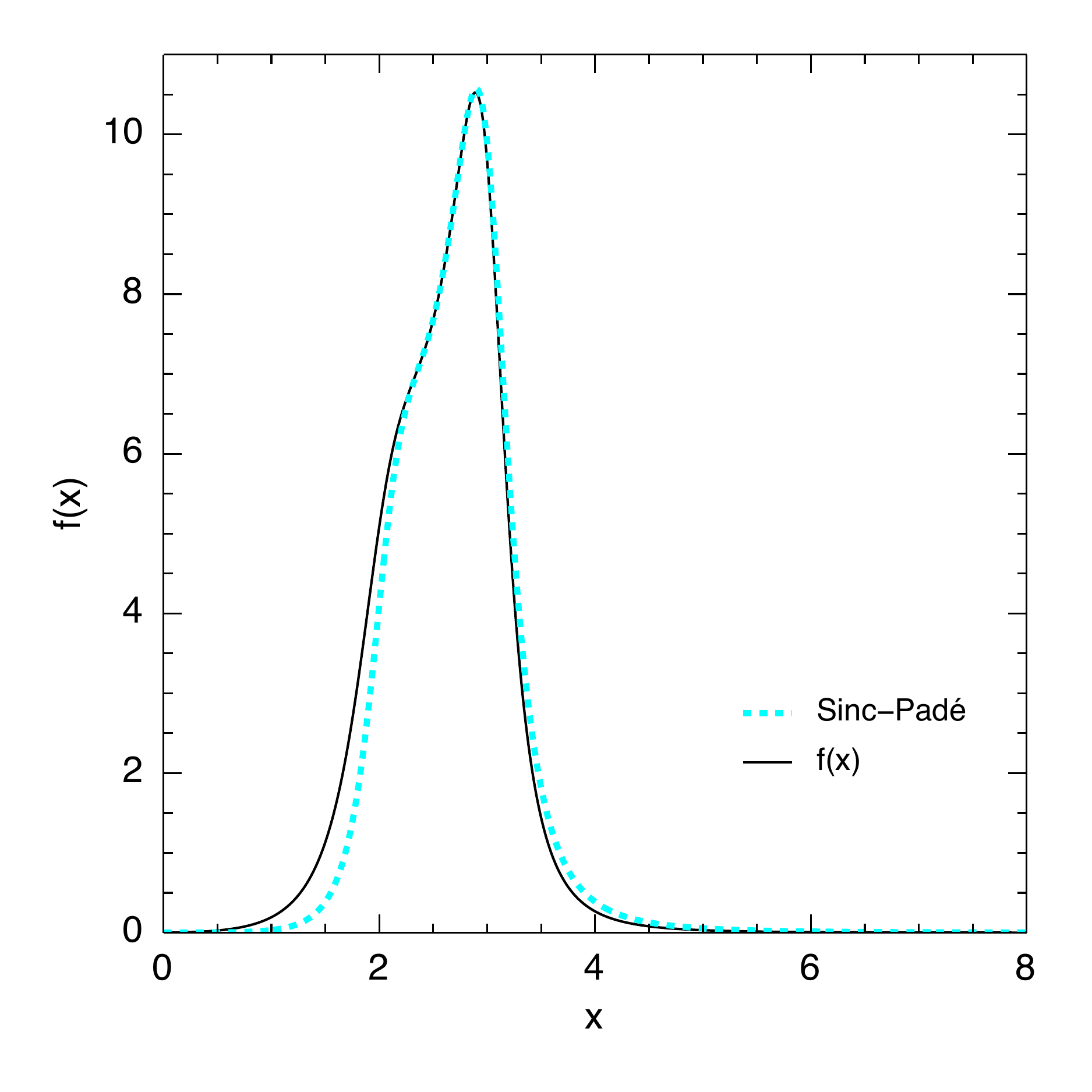}&
\includegraphics[width=0.45\textwidth]{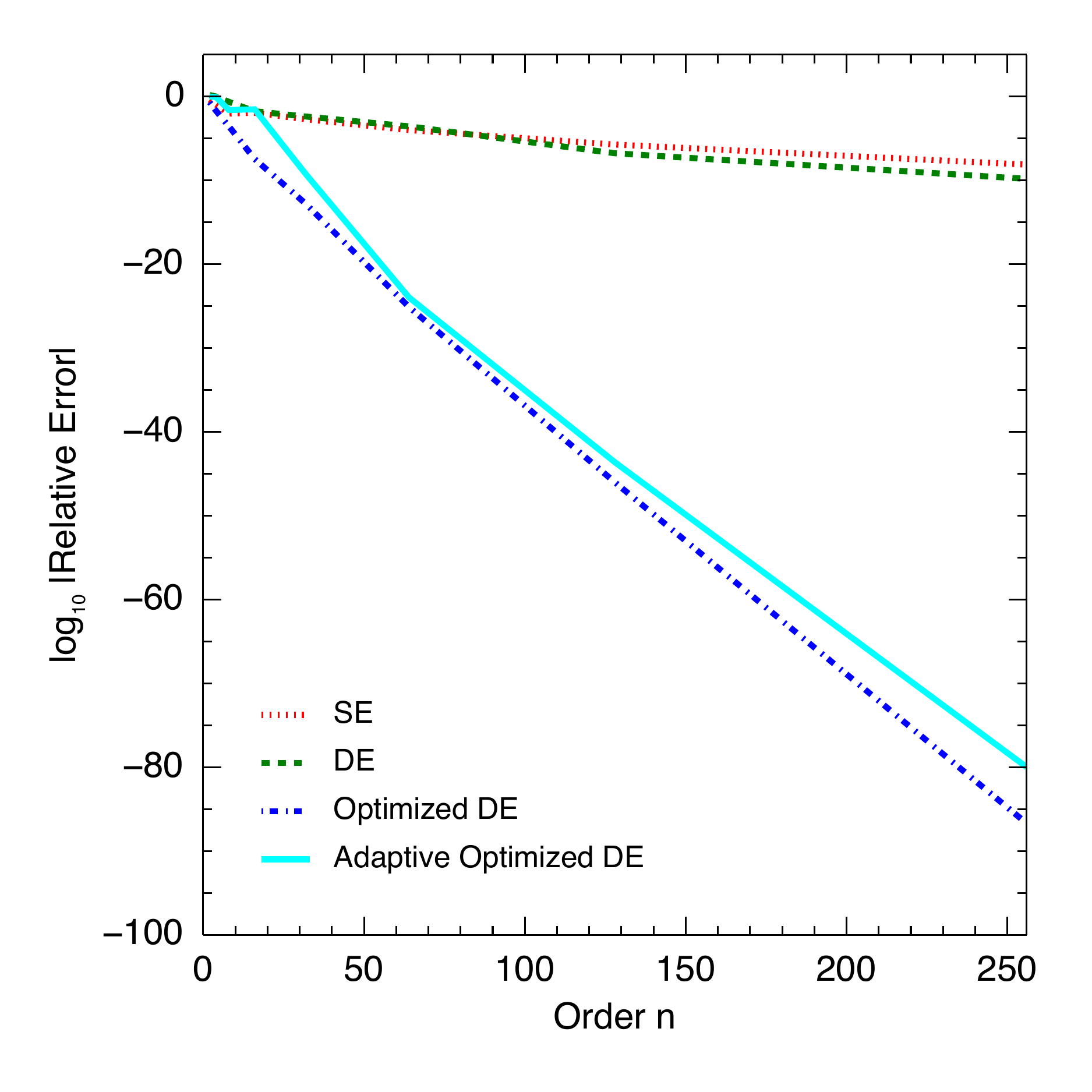}\\
(a) & (b)\\
\end{tabular}
\caption{In (a) the plot of the integrand of~\eqref{eq:SEandDEEx4} and the Sinc-Pad\'e approximant for $n=2^8$ and in (b) the performance of the trapezoidal rule with single, double, optimized double, and adaptive optimized double exponential variable transformations.}
\label{fig:Example4ploterr}
\end{center}
\end{figure}

In Figure~\ref{fig:Example4ploterr} (a), the integrand of~\eqref{eq:SEandDEEx4} is shown, and in Figure~\ref{fig:Example4ploterr} (b), the logarithm of the relative errors of the trapezoidal rule of order $n$ with single, double, and optimized double exponential variable transformations are plotted. The increase in convergence rate using the optimized variable transformation is a significant increase in efficiency over the double exponential transformation.

\section{Applications}

\subsection{Nonlinear waves}

For internal waves in stratified fluids of great depth, the Benjamin-Ono equation~\cite{Benjamin-29-559-67,Ono-39-1082-75} is a nonlinear partial integro-differential equation involving the Hilbert transform. While the KdV equation has soliton solutions that decay exponentially on the real line, the soliton solutions to the homogeneous Benjamin-Ono equation decay algebraically. The Hilbert transform on the real line is defined as~\cite{King-1-09}:
\begin{equation}\label{eq:HilbertTransform}
{\cal H}y(x) = \dfrac{1}{\pi}\dashint_{-\infty}^{+\infty} \dfrac{y(s)}{s-x}{\rm\,d}s,
\end{equation}
where the dash in the integral sign denotes the Cauchy principal value. The forced Benjamin-Ono equation can then be written as:
\begin{equation}\label{eq:BenjaminOno}
y_t + y y_x + {\cal H} y_{xx} = g(x-c\,t),\qquad x\in\mathbb{R},\quad t\ge0,\quad \lim_{x\to\pm\infty}y(x,t) = 0.
\end{equation}
for some $g$ with wave speed $c\in\mathbb{R}$. If we let $y(x,t) = y(x-c\,t)$, then we may consider the traveling wave solutions. Inserting such a substitution into~\eqref{eq:BenjaminOno} and integrating, we obtain:
\begin{align}
-c\,y' + y y' + {\cal H} y'' & = g(x-c\,t),\\
-c\,y + \dfrac{y^2}{2} + {\cal H} y' & = f(x-c\,t) = \int_{-\infty}^{x-c\,t}g(s){\rm\,d}s.
\end{align}
The traveling wave is then obtained by solving this equation for $t=0$.

The forced solutions to this equation have many properties that can be deduced from $f$. For example, if $f$ decays algebraically on the real line, the method of dominant balance shows that:
\begin{equation}
y(x) \sim -c^{-1}f(x),\quad{\rm as}\quad x\to\pm\infty,
\end{equation}
and will therefore behave similarly. As well, complex singularities in $f$ can potentially be found in $y$.

To continue, it is clear we require an approximation for the Hilbert transform. In~\cite{Stenger-93,Stenger-121-379-00}, Stenger derives a Sinc-based approximation for the Hilbert transform. We closely follow his development which is for a general variable transformation, and show how the optimized conformal map improves the convergence rate.

By making the invertible variable transformation $s=\phi(t):\mathbb{R}\to\mathbb{R}$ in~\eqref{eq:HilbertTransform}, we obtain:
\begin{equation}
{\cal H}y(x) = \dfrac{1}{\pi}\dashint_{-\infty}^{+\infty} \dfrac{y(\phi(t))\phi'(t)}{\phi(t)-x}{\rm\,d}t.
\end{equation}
The integrand multiplied by $t-\phi^{-1}(x)$ has but a removable singularity at $t = \phi^{-1}(x)$. Therefore, we may approximate with a Sinc basis:
\begin{equation}\label{eq:SincHilb}
\dfrac{y(\phi(t))\phi'(t)}{\phi(t)-x}(t-\phi^{-1}(x)) \approx \sum_{j=-n}^{+n} \dfrac{y(\phi(jh))\phi'(jh)}{\phi(jh)-x}(jh-\phi^{-1}(x))S(j,h)(t).
\end{equation}
Using the Hilbert transform~\cite{King-1-09}:
\begin{equation}
{\cal H} \dfrac{\sin x}{x} = \dfrac{1}{\pi}\dashint_{-\infty}^{+\infty} \dfrac{\sin s}{s(s-x)}{\rm\,d}s = \dfrac{\cos x - 1}{x},
\end{equation}
and dividing by the linear factor $t-\phi^{-1}(x)$ and integrating each basis function, the approximation~\eqref{eq:SincHilb} leads directly to:
\begin{equation}\label{eq:SincHilbert}
{\cal H}y(x) \approx \dfrac{h}{\pi}\sum_{j=-n}^{+n}y(\phi(jh))\phi'(jh) \dfrac{\cos\left(\pi (\frac{\phi^{-1}(x)}{h}-j)\right)-1}{x-\phi(jh)}.
\end{equation}

Using the Sinc-based approximation for the Hilbert transform, the solution of the forced Benjamin-Ono equation can be obtained by approximating the solution $y(x)$ in a Sinc basis:
\begin{equation}
y(x) \approx y_n(x;\phi) = \sum_{j=-n}^{+n} y_j S(j,h)(\phi^{-1}(x)),
\end{equation}
and solving the system of nonlinear equations obtained by collocating at the Sinc points $x_k = \phi(kh)$:
\begin{equation}
-c\,y_n(k\,h;\phi) + y_n^2(k\,h,\phi) + {\cal H}y_n'(k\,h,\phi) = f(k\,h),\qquad k=-n,\ldots,n,
\end{equation}
for the $2n+1$ unknowns $\{y_j\}_{|j|\le n}$ by Newton iteration.

For the purposes of illustration, we consider the functions $f(x)$ which yield the solutions:
\begin{equation}
y(x) = \sum_{i=1}^m \dfrac{\epsilon_i^2}{((x-\delta_i)^2+\epsilon_i^2)},
\end{equation}
for the values $m=3$, $\delta_1+{\rm i}\epsilon_1 = -1+0.3{\rm i}$, $\delta_2+{\rm i}\epsilon_2 = 0+0.1{\rm i}$, and $\delta_3+{\rm i}\epsilon_3 = 1+0.2{\rm i}$. This allows for an exact comparison with an analytic expression. In addition, the wave speed $c=1$ is used. Table~\ref{table:BenjaminOno} summarizes the variable transformations and the parameters used.

\begin{table}[htbp]
\begin{center}
\caption{Transformations and parameters for~\eqref{eq:BenjaminOno}.}
\label{table:BenjaminOno}
\begin{tabular*}{\hsize}{@{\extracolsep{\fill}}c|ccc}
\hline
& Single & Double & Optimized Double \\
\hline
$\phi(t)$ & $\sinh(t)$ & $\sinh(\pi/2\sinh(t))$ & $\sinh(h(t))$\\
$\rho$ or $\gamma$ & $1$ & $1$ & $1$\\
$\beta$ or $\beta_2$ & $1/2$ & $\pi/4$ & $5.7257\times10^{-8}$\\
$d$ & $0.10017$ & $0.06381$ & $\pi/2$\\
\hline
\end{tabular*}
\end{center}
\end{table}

In addition, the optimized transformation is given by:
\begin{equation}
h(t) \approx 1.1451\times10^{-7}\sinh(t) + 0.04531 + 0.06359\,t - 1.2134\times10^{-4}t^2.
\end{equation}

\begin{figure}[htbp]
\begin{center}
\begin{tabular}{cc}
\includegraphics[width=0.45\textwidth]{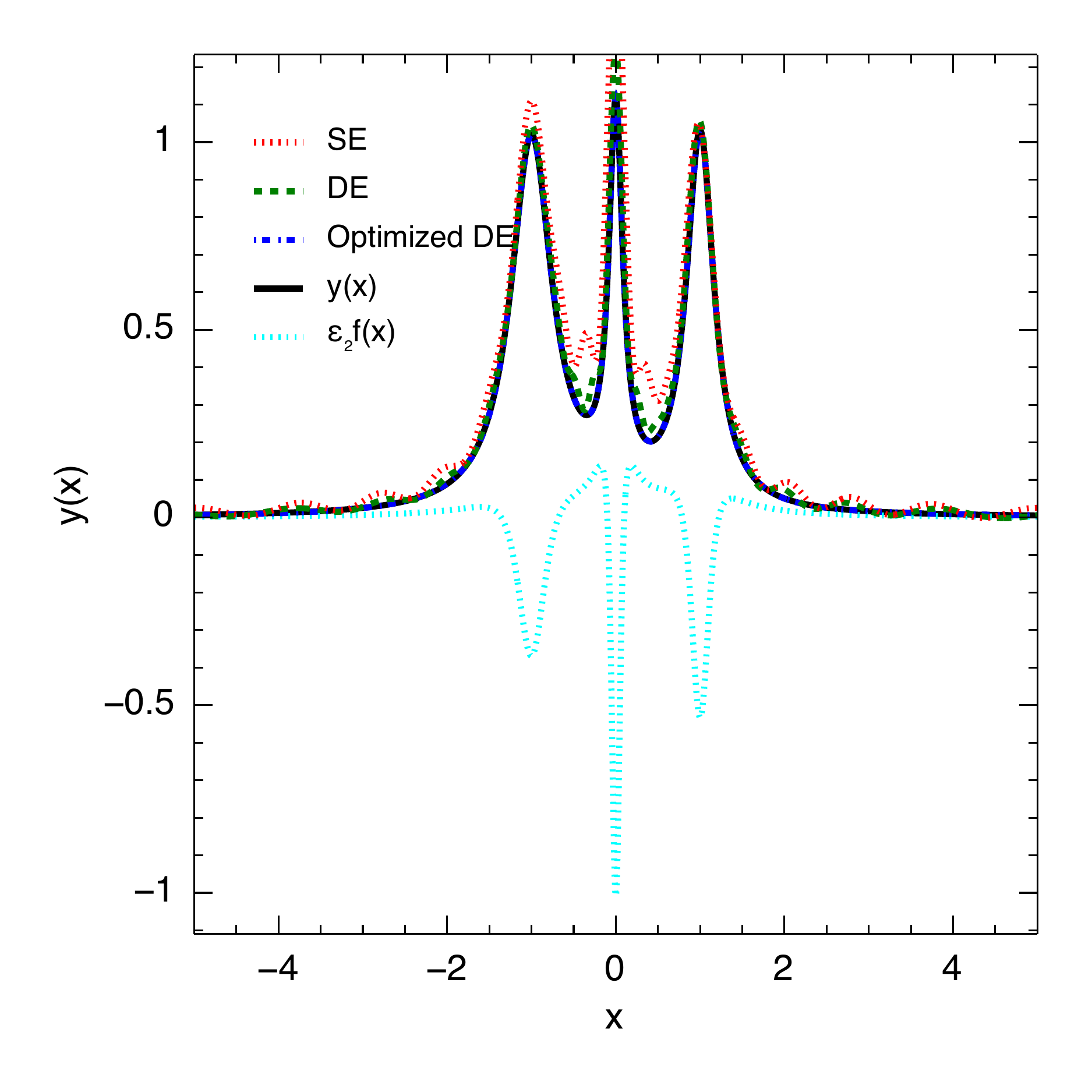}&
\includegraphics[width=0.45\textwidth]{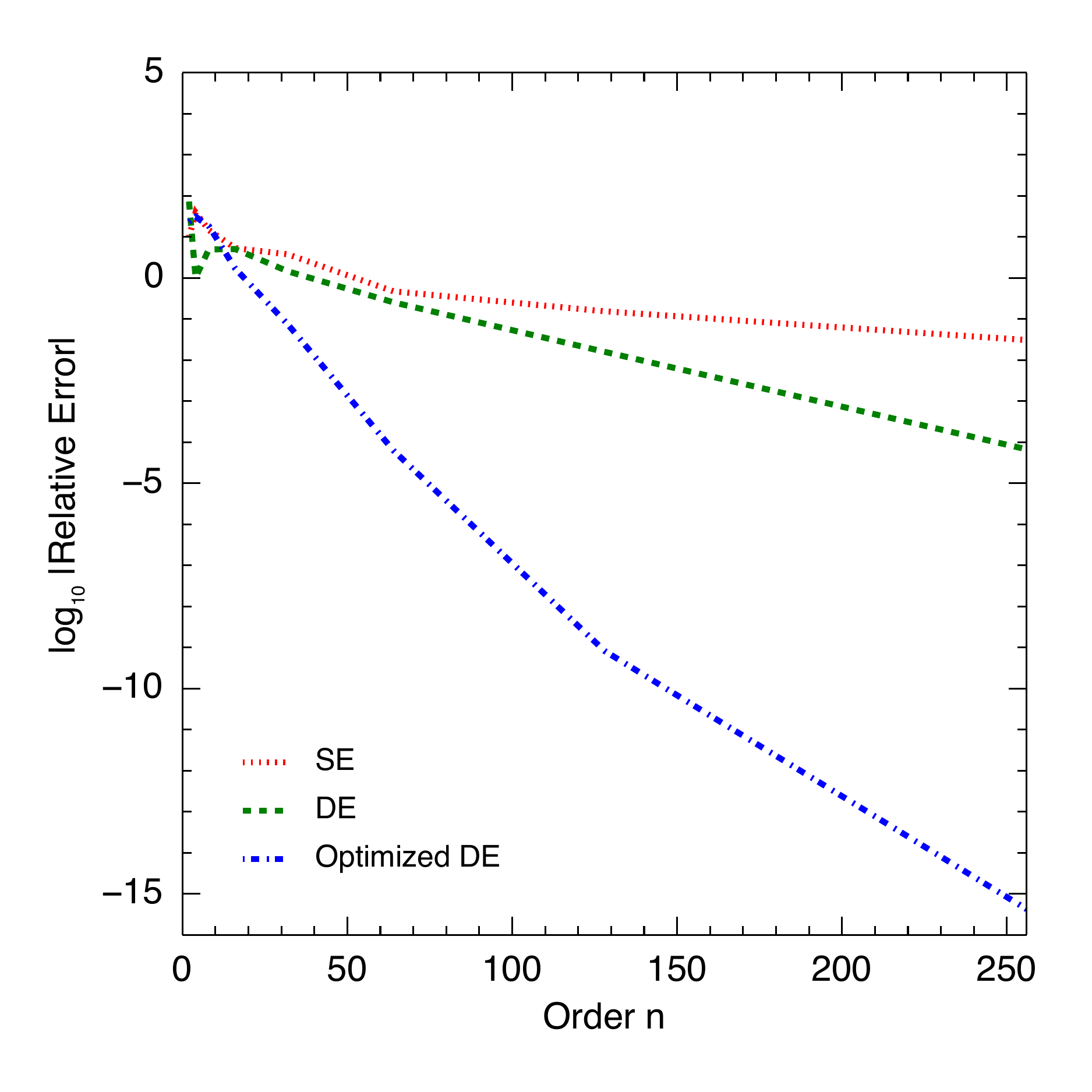}\\
(a) & (b)\\
\end{tabular}
\caption{In (a) the solution of~\eqref{eq:BenjaminOno} along with the three Sinc approximations for $n=2^{5}$ and the scaled forcing function and in (b) the performance of the Sinc approximation with single, double, and optimized double exponential variable transformations.}
\label{fig:BenjaminOno}
\end{center}
\end{figure}

In Figure~\ref{fig:BenjaminOno}, we approximate the relative error:
\begin{equation}
\sup_{x\in\mathbb{R}}\left| \dfrac{y(x)-y_n(x;\phi)}{y(x)}\right|,
\end{equation}
by computing the maximum of the difference and quotient at $101$ equally spaced points in the interval $[-5,5]$. The inverse optimized map $\phi^{-1}$ is conveniently computed via Newton iteration, as the map and its first derivative are already required in the system of collocated equations. The increase in convergence rate using the optimized variable transformation is a significant increase in efficiency over the double exponential transformation.

\subsection{Multidimensional integrals}

There are many applications in physics that warrant the evaluation of $m$-dimensional integrals. Examples we are interested in include: magnetic susceptibility integrals in the Ising theory of solid-state physics~\cite{Bailey-Borwein-Crandall-39-12271-06}, which form terms in a series that represents the dependence of magnetic susceptibility on temperature; and, box integrals~\cite{Bailey-Borwein-Crandall-206-196-07}, which are essentially $m$-dimensional expectations of the $s^{\rm th}$ power of the distance in a unit hypercube $\langle |\vec{r}|^{s}\rangle_{\vec{r}\in[0,1]^m}$. Such integrals have applications in probability theory, and in potential theory such as ``jellium'' potentials. After making substantial analytic advances in the theory of box integrals in~\cite{Bailey-Borwein-Crandall-79-1839-10}, the authors acknowledge that some of the analytical techniques they used would not apply to other $m$-dimensional expectations, and posit that $\langle e^{-\kappa|\vec{r}|}\rangle_{\vec{r}\in[0,1]^m}$ should remain extremely difficult to evaluate in any general way. In this section, we construct such a general way to calculate these integrals. Explicitly:
\begin{equation}
\langle e^{-\kappa|\vec{r}|}\rangle_{\vec{r}\in[0,1]^m} = \int_{[0,1]^m} e^{-\kappa (r_1^2+\cdots+r_m^2)^{1/2}}{\rm\,d}r_1\cdots{\rm\,d}r_m.\label{eq:BoxExp}
\end{equation}
Using the same dimensional-reduction technique in~\cite{Bailey-Borwein-Crandall-206-196-07}, and the incomplete gamma function~\cite[\S 8.350 1.]{Gradshteyn-Ryzhik-07}:
\begin{equation}
\gamma(m,a) = a^m\int_0^1 x^{m-1}e^{-a\,x}{\rm\,d}x = (m-1)! - e^{-a}\sum_{j=0}^{m-1}\binom{m-1}{j}\dfrac{j!}{a^{j-m+1}},
\end{equation}
we obtain:
\begin{align}
\langle e^{-\kappa|\vec{r}|}\rangle_{\vec{r}\in[0,1]^m} = \dfrac{m}{2^{m-1}}&\int_{[-1,1]^{m-1}}{\rm\,d}x_1\cdots{\rm\,d}x_{m-1}\left(\dfrac{(m-1)!}{\kappa^m(x_1^2+\cdots+x_{m-1}^2+1)^{m/2}}\right.\nonumber\\
&\left. - \sum_{j=0}^{m-1}\binom{m-1}{j}\dfrac{j!e^{-\kappa(x_1^2+\cdots+x_{m-1}^2+1)^{1/2}}}{\kappa^{j+1}(x_1^2+\cdots+x_{m-1}^2+1)^{(j+1)/2}}\right).\label{eq:BoxExpRedm1}
\end{align}
multidimensional integrals are a challenge for univariate numerical integration methods because of the curse of dimensionality, whereby the dimension $m$ increases the number of sample points of the one-dimensional case $N=2n+1$ geometrically as ${\cal O}(N^m)$, reaching the limits of modern computational power quite quickly. Nevertheless, positive results may be obtained for lower dimensions, especially for extremely high accuracy, for which even tens of digits qualifies in this setting. We compare and contrast the trapezoidal rule on~\eqref{eq:BoxExp} for the single, double, and optimized double exponential transformations. These transformations are summarized in Table~\ref{table:BoxExp}. Of course, $m-1$ transformations $\{\phi_{\ell}(t_1,\ldots,t_{m-1})\}_{\ell=1}^{m-1}$ will be required to induce decay at all the boundaries of the hypercube in~\eqref{eq:BoxExp}.

\begin{table}[htbp]
\begin{center}
\caption{Transformations and parameters for~\eqref{eq:BoxExpRedm1}.}
\label{table:BoxExp}
\begin{tabular*}{\hsize}{@{\extracolsep{\fill}}c|ccc}
\sphline
& Single & Double & Optimized Double \\
\sphline
$\phi_{\ell}(t_1,\ldots,t_{m-1})$ & $\tanh(t_{\ell}/2)$ & $\tanh(\pi/2\sinh(t_{\ell}))$ & {\footnotesize$\tanh\left(\tan^{-1}\left(\sqrt{\sum_{j=1}^{\ell-1}\phi_j^2+1}\right)\sinh(t_{\ell})\right)$}\\
$\rho$ or $\gamma$ & $1$ & $1$ & $1$\\
$\beta$ or $\beta_2$ & $1$ & $\pi/2$ & $\pi/4$\\
$d$ & $\pi/2$ & $\pi/6$ & $\pi/2$\\
\hline
\end{tabular*}
\end{center}
\end{table}

Note that while the $\ell^{\rm th}$ optimized double exponential transformation calls all previous ones through the term $\sqrt{\sum_{j=1}^{\ell-1}\phi_j^2+1}$, this comes at no extra cost because it is extracted during the construction of the term $x_1^2+\cdots+x_{m-1}^2+1$, which occurs in~\eqref{eq:BoxExpRedm1}.

\begin{figure}[htbp]
\begin{center}
\begin{tabular}{ccc}
\includegraphics[width=0.3\textwidth]{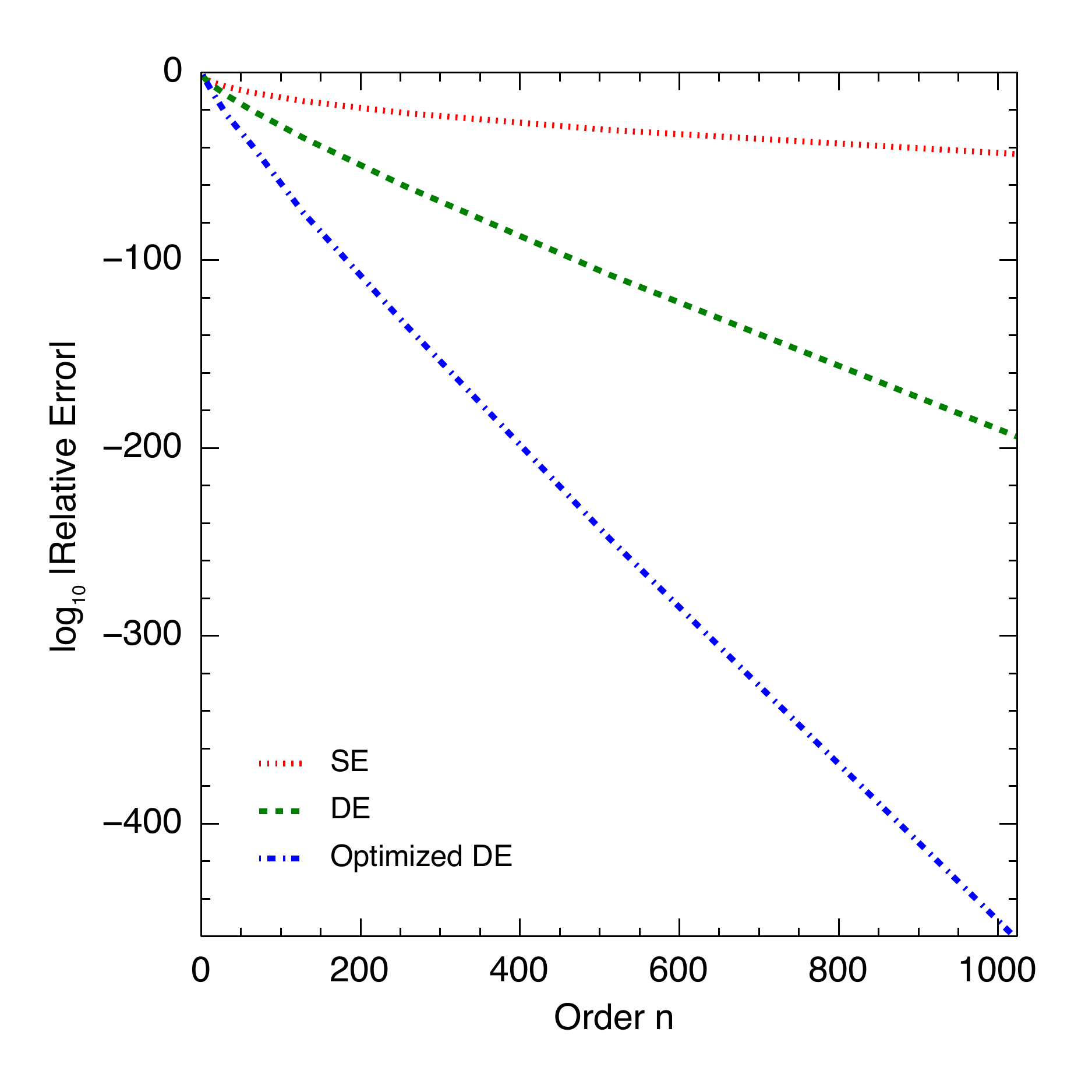}&
\includegraphics[width=0.3\textwidth]{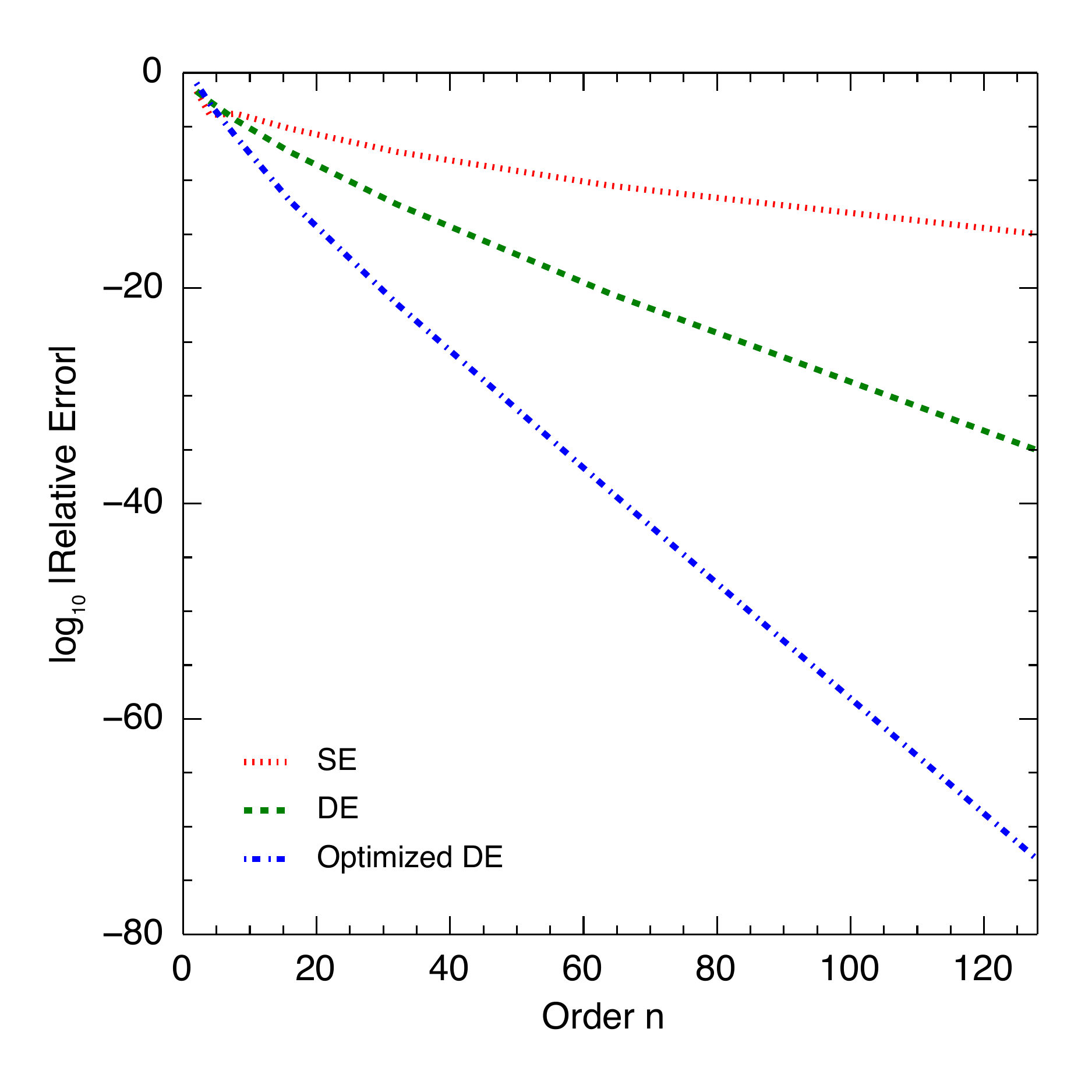}&
\includegraphics[width=0.3\textwidth]{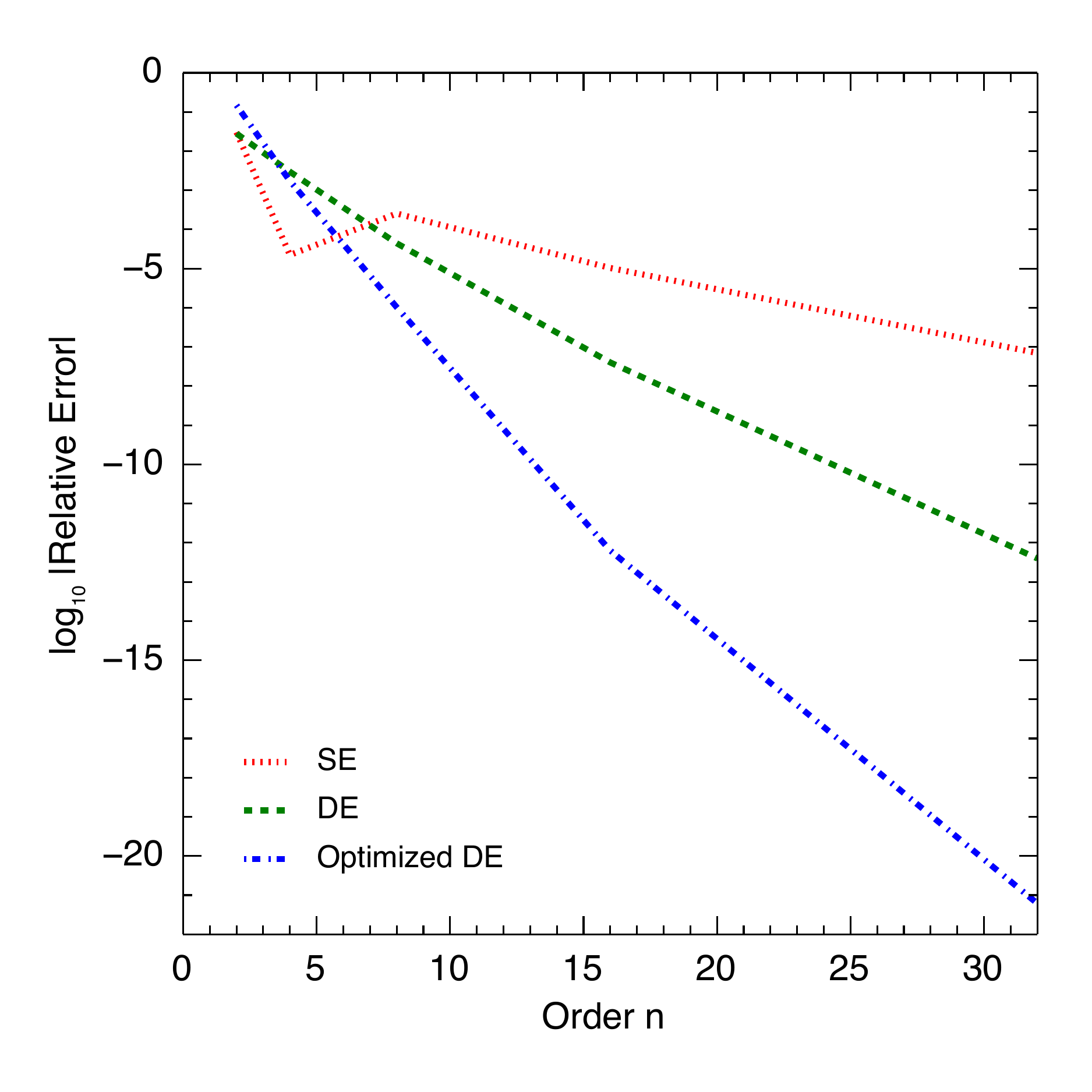}\\
(a) & (b) & (c)\\
\end{tabular}
\caption{The performance of the trapezoidal rule with single, double, and optimized double exponential variable transformations. In all figures, $\kappa=1$, and in (a) $m=2$, in (b) $m=3$, and in (c) $m=4$.}
\label{fig:BoxExperr}
\end{center}
\end{figure}

In Figure~\ref{fig:BoxExperr}, the logarithm of the relative errors of the trapezoidal rule of order $n$ with single, double, and optimized double exponential variable transformations are plotted. The increase in convergence rate using the optimized variable transformation is a significant increase in efficiency over the double exponential transformation.

The dimension of the so-called box integrals of~\cite{Bailey-Borwein-Crandall-206-196-07} was reduced completely to one in every case by using an integral representation of the integrand $|\vec{r}|^s$, which allowed the integrals over $r_1$, $r_2$, etc\ldots~to be separated and written as the $m^{\rm th}$ power of the error function. The authors postulated the exponential expectation to be a challenge because no integral representation was known to reduce the dimension of the problem. However, using the Bessel representation for $K_{-1/2}(\kappa|\vec{r}|)$ from~\cite[\S 8.432 7.]{Gradshteyn-Ryzhik-07}:
\begin{equation}
e^{-\kappa|\vec{r}|} = \sqrt{\dfrac{2\,\kappa|\vec{r}|}{\pi}}K_{-1/2}(\kappa|\vec{r}|) = \sqrt{\dfrac{\kappa}{2\pi}}\int_0^\infty \dfrac{e^{-\frac{\kappa}{2}(t+\frac{r^2}{t})}}{t^{1/2}}{\rm\,d}t,
\end{equation}
and the error function representation~\cite[\S 3.321 1.]{Gradshteyn-Ryzhik-07}:
\begin{equation}
{\rm erf}(u) = \dfrac{2\,u}{\sqrt{\pi}}\int_0^1e^{-u^2x^2}{\rm\,d}x,
\end{equation}
we are able to obtain the formula:
\begin{equation}
\langle e^{-\kappa|\vec{r}|}\rangle_{\vec{r}\in[0,1]^m} = \dfrac{1}{2}\left(\dfrac{\pi}{2\,\kappa}\right)^{\frac{m-1}{2}}\int_0^\infty t^{\frac{m-1}{2}}e^{-\kappa\,t/2}{\rm erf}^m\left(\sqrt{\dfrac{\kappa}{2\,t}}\right){\rm\,d}t.\label{eq:BoxExpRed1}
\end{equation}
Because these are such challenging integrals, we include rounded approximate values in Table~\ref{table:BoxExpVals} for the sake of reproducibility. The one-dimensional integral~\eqref{eq:BoxExpRed1} allows for the calculation of the high precision results, and it is also used for the calculation of the relative error in Figure~\ref{fig:BoxExperr}.

\begin{table}[htbp]
\begin{center}
\caption{Numerical Evaluation of~\eqref{eq:BoxExp} using~\eqref{eq:BoxExpRed1}.}
\label{table:BoxExpVals}
\begin{tabular}{cc|l}
\sphline
$m$ & $\kappa$ & \multicolumn{1}{c}{$\langle e^{-\kappa|\vec{r}|}\rangle_{\vec{r}\in[0,1]^m}$} \\
\sphline
$2$ & $1.0$ & 0.48499~93872~72994~84128~76561~86058~31858~19718\\
$3$ & $1.0$ & 0.39822~04526~88323~04659~07885~63033~98432~76981\\
$4$ & $1.0$ & 0.33843~80876~94843~90404~45300~56568~55958~16022\\
$5$ & $1.0$ & 0.29379~80818~76007~61424~12657~48176~65958~00955\\
\hline
\end{tabular}
\end{center}
\end{table}

\section{Numerical Discussion}

The numerical experiments are all programmed in Julia~\cite{Julia-12}, calling at times GNU's MPFR library for arbitrary precision arithmetic, OpenBLAS for solving the linear systems and Ipopt~\cite{Wachter-Biegler-106-25-06} for solving the nonlinear program~\eqref{eq:SCaltparamprob}. As can be seen in the figures showing the relative errors, the maximization of the convergence rates provides a significant improvement for the double exponential variable transformations. With an equal number of function evaluations, the optimized double exponential formulas provide approximately $2.5$--$4$ times as many correct digits. The conformal maps achieve this by locating singularity pre-images on the boundary of the strip $\partial\mathscr{D}_{\frac{\pi}{2}}$. 

In Examples~\ref{example:SEandDEEx1}--\ref{example:SEandDEEx4}, one integral is treated on each of the canonical domains: in Example~\ref{example:SEandDEEx1}, two different endpoint and two pairs of different near-contour singularities are treated; in Example~\ref{example:SEandDEEx2}, four pairs of different near-contour singularities are treated; in Example~\ref{example:SEandDEEx3}, an infinite array of singularities is treated; and, in Example~\ref{example:SEandDEEx4}, the adaptive optimized method is shown to successfully approximate the loci of the three pairs of near-contour singularities.

In every case, the nonlinear program~\eqref{eq:SCaltparamprob} still does not ensure analyticity in the strip $\mathscr{D}_{\frac{\pi}{2}}$. In Examples~\ref{example:SEandDEEx1},~\ref{example:SEandDEEx2} and~\ref{example:SEandDEEx4}, the compositions $\psi(h(t))$ actually cross the branches of the square root functions, and in Example~\ref{example:SEandDEEx3}, there are even more poles than are shown in Figure~\ref{fig:Example3maps}. Yet still, a significant increase in convergence is observed. It is quite clear that some singular effects are numerically more damaging than others.

Example~\ref{example:SEandDEEx4} shows the use of the Sinc-Pad\'e approximants for the adaptive optimized algorithm~\ref{alg:adaptiveDE}. While in Figure~\ref{fig:Example4ploterr}, the Sinc-Pad\'e approximants serve as relatively poor approximations of the integrand, their ability to approximate singular points is acceptable. In Table~\ref{table:SEandDEEx4SincPade}, the Sinc-Pad\'e approximants obtain $2$-$3$ correct digits for the first order poles $\delta_2\pm{\rm i}\epsilon_2$ and $\delta_3\pm{\rm i}\epsilon_3$, but struggle to obtaining an accurate location of the branch points $\delta_1\pm{\rm i}\epsilon_1$. This is entirely related to the rational limitations of the Sinc-Pad\'e approximants, and suggests that approximating essential singularities, for example, surpasses the capabilities of the Pad\'e methods.

While the nonlinear program~\eqref{eq:SCaltparamprob} is successful in the current endeavours, further research in other conformal maps would be fruitful. For example, it is still unclear how to maximize the convergence rate in cases of countably infinite singularities, such as:
\begin{equation}
\int_{-\infty}^{+\infty} \dfrac{\tanh x}{x(1+x^2)}{\rm\,d}x.
\end{equation}
A single exponential transformation such as $\frac{\pi}{2}\sinh$ creates an array of singularities along $\pm\frac{{\rm i}\pi}{2}$. Therefore, any further composition will inevitably cause the limits of this array to approach the real axis without bound. The problem has been discussed in~\cite{Okayama-Tanaka-Matsuo-Sugihara-125-511-13,Tanaka-Okayama-Matsuo-Sugihara-125-545-13}, and the result is almost the same convergence property as a single exponential transformation.

Nevertheless, the applications show how the conformal maps can be used to obtain substantial increases in accuracy. For the forced Benjamin-Ono equation for nonlinear waves, only the optimized double exponential transformation is able to provide any accuracy that resembles the solution even for a large basis. And for the multidimensional integrals, the gain in accuracy is also significant.

\bibliography{/Users/Mikael/Mik}

\end{document}